\RequirePackage{ifthen}
\newboolean{MPA}
\setboolean{MPA}{false}

\ifthenelse {\boolean{MPA}}
{
\documentclass{svjour3}
\smartqed
} {
\documentclass[11pt]{article}
\usepackage[margin=1.5in]{geometry}
}

\usepackage{graphicx}
\usepackage{amsmath}
\usepackage{amssymb}
\usepackage{hyperref}
\usepackage{color}

\usepackage{mathrsfs} 

\usepackage{xspace} 

\usepackage{soul} 

\usepackage[makeroom]{cancel} 

\usepackage{algorithm}
\usepackage[noend]{algpseudocode}

\floatname{algorithm}{}

\newcommand{\R}{\mathbb R}
\newcommand{\Z}{\mathbb Z}
\newcommand{\Q}{\mathbb Q}

\newcommand{\id}{I}

\newcommand{\doubletilde}[1]{\tilde{\raisebox{0pt}[0.85\height]{$\tilde{#1}$}}}

\let\st\relax
\DeclareMathOperator{\st}{s.t.}

\DeclareMathOperator{\proj}{proj}
\DeclareMathOperator{\spn}{span}

\DeclareMathOperator{\width}{width}

\newcommand{\instances}{\mathscr I}
\newcommand{\solutions}{\mathscr A}

\newcommand{\ie}{i.e., }

\renewcommand{\P}{\mathcal P}
\newcommand{\B}{\mathcal B}
\newcommand{\C}{\mathcal C}
\renewcommand{\L}{\mathcal L}
\newcommand{\M}{\mathcal M}
\newcommand{\N}{\mathcal N}
\renewcommand{\S}{\mathcal S}

\newcommand{\E}{\mathcal E}

\newcommand{\G}{G}
\renewcommand{\H}{H}
\newcommand{\zcoeff}{l}
\newcommand{\lastrowR}{\hat r}
\newcommand{\lastrowU}{u}
\newcommand{\direction}{v}
\newcommand{\did}{d}
\newcommand{\basmatB}{B}

\newcommand{\floor}[1]{\lfloor#1\rfloor}

\newcommand{\round}[1]{\lfloor#1\rceil}

\newcommand{\ceil}[1]{\lceil#1\rceil}
\newcommand{\ceilL}[1]{\left\lceil#1\right\rceil}

\newcommand{\abs}[1]{\lvert#1\rvert}
\newcommand{\absL}[1]{\left\lvert#1\right\rvert}

\newcommand{\norm}[1]{\lVert#1\rVert}
\newcommand{\normL}[1]{\left\lVert#1\right\rVert}

\newcommand{\Fnorm}[1]{\lVert#1\rVert_F}

\newcommand{\pare}[1]{\left(#1\right)}
\newcommand{\bra}[1]{\left\{#1\right\}}

\newcommand{\transp}{\mathsf T}

\newcommand{\ctr}{a}
\newcommand{\low}{\ell}
\newcommand{\upp}{u}

\newcommand{\constlen}{2^{p(p-1)/4}}
\newcommand{\constdec}{q}
\newcommand{\constelli}{r}
\newcommand{\constwidth}{s}
\newcommand{\constboxes}{\varphi}

\newcommand{\linearized}{g}
\newcommand{\linearizedinlemma}{q'}

\usepackage{comment}
\specialcomment{knownproof}{\color{blue}}{\color{black}}
\excludecomment{knownproof}

\newcounter{step}
\newcommand{\step}[1]{\refstepcounter{step} \smallskip \noindent \textbf{Step~\thestep: #1.}}

\ifthenelse {\boolean{MPA}}
{

\newenvironment{prf}[1][]
{\begin{proof}}
{\qed \end{proof}}

\newenvironment{prfh}[1][]
{\begin{proof}}
{\end{proof}}

\newcommand{\qedhere}{\tag*{\qed}}

\newcounter{claim} 
\renewenvironment{claim}[1][]
{\refstepcounter{claim} \begin{trivlist} \item[] {\bf Claim~\theclaim}\space#1 \itshape}
{\end{trivlist}}

\journalname{Mathematical Programming A}

} {

\usepackage{amsthm}
\newtheorem{theorem}{Theorem}
\newtheorem{proposition}{Proposition}
\newtheorem{lemma}{Lemma}
\newtheorem{corollary}{Corollary}

\newtheorem{claim}{Claim}

\newenvironment{prf}[1][]
{\begin{proof}}
{\end{proof}}

\newenvironment{prfh}[1][]
{\begin{proof}}
{\end{proof}}

}

\newtheorem{observation}{Observation}

\begin{document}

\title{An Approximation Algorithm for Indefinite \\ Mixed Integer Quadratic Programming}

\ifthenelse {\boolean{MPA}}
{
\titlerunning{An Approximation Algorithm for Indefinite MIQP}

\author{Alberto Del Pia}
\institute{Alberto~Del~Pia \at
              Department of Industrial and Systems Engineering 
              \& Wisconsin Institute for Discovery \\
              University of Wisconsin-Madison, Madison, WI, USA \\
              \email{delpia@wisc.edu}}
}
{
\author{Alberto Del Pia
\thanks{Department of Industrial and Systems Engineering \& Wisconsin Institute for Discovery,
             University of Wisconsin-Madison, Madison, WI, USA.
             E-mail: {\tt delpia@wisc.edu}.}}
}

\date{\today}

\maketitle

\begin{abstract}
In this paper, we give an algorithm that finds an $\epsilon$-approximate solution to a mixed integer quadratic programming (MIQP) problem.
The algorithm runs in polynomial time if the rank of the quadratic function and the number of integer variables are fixed.
The running time of the algorithm is expected unless P=NP.
In order to design this algorithm we introduce the novel concepts of spherical form MIQP and of aligned vectors, and we provide a number of results of independent interest.
In particular, we give a strongly polynomial algorithm to find a symmetric decomposition of a matrix, and show a related result on simultaneous diagonalization of matrices.
\ifthenelse {\boolean{MPA}}
{
} {}
\end{abstract}

\ifthenelse {\boolean{MPA}}
{}{
}



\section{Introduction}


\emph{Mixed Integer Quadratic Programming} (MIQP) is an optimization problem where the objective function is a general quadratic function, the constraints are linear inequalities, and some of the variables are required to be integers.
Formally, given a symmetric matrix $H \in \Q^{n \times n}$, a matrix $W \in \Q^{m \times n}$, vectors $h \in \Q^n$, $w \in \Q^m$, and $p \in \{0,1,\dots,n\}$, we seek a vector $x \in \R^n$ that attains
\begin{align}
\label{prob main}
\tag{MIQP}
\begin{split}
\min & \quad x^\transp H x + h^\transp x \\
\st & \quad Wx \le w \\
& \quad x \in \Z^p \times \R^{n-p}. \\
\end{split}
\end{align}

Many important applications can be modeled as MIQPs,
in areas such as operations research, engineering, computer science, physics, biology, finance, economics, and artificial intelligence.
MIQP reduces to \emph{Mixed Integer Linear Programming} (MILP) when $H$ is a zero matrix, and to \emph{Quadratic Programming} (QP) if $p=0$.
Moreover, MIQP is a prototypical \emph{Mixed Integer Nonlinear Programming} (MINLP) problem, as it captures the critical elements of those models,
but in the simplest possible way, making it the natural first step to construct efficient algorithms for MINLP.

The decision version of MIQP lies in the complexity class NP~\cite{dPDeyMol17}.
Furthermore, 
MIQP is strongly NP-hard~\cite{GarJohSto76}, and remains NP-hard even if $H$ has rank one and $p=0$~\cite{ParVav91}. 
This implies the lack of efficient algorithms for solving this class of optimization problems in its full generality. 
%
%

The main result of this paper is an approximation algorithm for MIQP.
In order to state our result, we first give the definition of~$\epsilon$-approximate solution.
%
%
%
%
%
%
%
Consider an instance of \eqref{prob main}, 
and assume 
that it has an optimal solution $x^*$.
Let $f(x)$ denote the objective function,
and let $f_{\max}$ be the maximum value of~$f(x)$ on the feasible region.
For $\epsilon \in [0,1]$, we say that a feasible point $x^\diamond$ is an \emph{$\epsilon$-approximate solution} if
\begin{equation*}
f(x^\diamond) - f(x^*) \le \epsilon \cdot (f_{\max} - f(x^*)).
\end{equation*}
Note that only optimal solutions are $0$-approximate solutions, while any feasible point is a $1$-approximate solution.
The definition of~$\epsilon$-approximate solution has some useful invariance properties which make it a natural choice in this setting. 
For instance, it is preserved under dilation and translation of the objective function, and it is insensitive to affine transformations of the objective function and of the feasible region, like for example changes of basis.
Our definition of approximation has been used in earlier works, and we refer to~\cite{NemYud83,Vav92c,BelRog95,KleLauPar06} for more details.
We can now state our main result.

\begin{theorem}
\label{th main}
For every $\epsilon \in (0,1]$, there is an algorithm that finds an $\epsilon$-approximate solution to a bounded \eqref{prob main}, if it exists. 
The running time of the algorithm is polynomial in the size of the input and in~$1/\epsilon$, provided that the rank $k$ of the matrix $H$ and the number of integer variables $p$ are fixed numbers.
\end{theorem}

This is the first known polynomial time approximation algorithm for MIQP with $k$ and $p$ fixed.
In particular, note that the dimension $n$ of the problem is not required to be fixed.
The running time of the algorithm exhibits a polynomial dependence on the size of the instance and on $1/\epsilon$, and an exponential dependence on $k$ and on $p$.
It is known that 
this dependence 
is expected unless P=NP, and we refer the reader to the discussion below the statement of Theorem~1 in~\cite{dP18}.



One might wonder if the boundedness assumption can be relaxed in Theorem~\ref{th main}, with the understanding that, if the input MIQP is unbounded, the algorithm should return at least a statement that the instance is unbounded.
The next theorem implies that the boundedness assumption cannot be removed unless P=NP.


\begin{theorem}
\label{th hard}
Determining whether \eqref{prob main} is unbounded is NP-complete, even if the rank $k$ of the matrix $H$ equals three and the number $p$ of integer variables is zero.
\end{theorem}

\begin{prf}
From Theorem~4 in~\cite{dPDeyMol17}, the decision problem in the statement is in NP, thus we only need to show the NP-hardness.
In Sections 2 and 3 in~\cite{ParVav91}, the authors present a QP of the form 
\begin{align}
\label{pr ParVav}
\min \{ x_1 - x_2^2 : Wx \le w, \ x \in \R^n\}
\end{align}
with nonnegative optimum objective value, and for which it is NP-hard to determine if the optimum value is zero.
Since every bounded QP has an optimal solution of polynomial size~\cite{Vav90}, there is a number $\phi$ which is polynomial in the size of the input QP \eqref{pr ParVav} for which the optimum objective value is either zero or strictly larger than $2^{-\phi}$.

Consider now the QP 
\begin{align}
\label{pr hard}
\min \{ x_1 x_{n+1} - x_2^2 - 2^{-\phi} x_{n+1}^2: ( W \mid -w ) x \le 0, \ x_{n+1} \ge 1, \ x \in \R^{n+1}\}.
\end{align}
Notice that the rank of the objective function is three. Thus, to conclude the proof of the theorem we only need to show that \eqref{pr hard} is unbounded if and only if the optimum value of \eqref{pr ParVav} is zero.

Assume that the optimum value of \eqref{pr ParVav} is zero.
Then there is a point $\bar x \in \R^n$ with $W \bar x \le w$ and $\bar x_1 - \bar x_2^2=0$.
Consider now the set of vectors in~$\R^{n+1}$ given by $(\lambda \bar x, \lambda)$, for $\lambda \ge 1$.
Note that all these vectors are feasible to \eqref{pr hard}.
Furthermore, the objective value of~$(\lambda \bar x, \lambda)$ is $\lambda^2 (\bar x_1 - \bar x_2^2 - 2^{-\phi}) = - \lambda^2 2^{-\phi}$ which goes to $-\infty$ as $\lambda \to \infty$.
Therefore, \eqref{pr hard} is unbounded.

Next, assume that the optimum value of \eqref{pr ParVav} is positive, therefore strictly larger than $2^{-\phi}$.
Consider a vector feasible to \eqref{pr hard} and note that it can be written as $(\bar \lambda \bar x, \bar \lambda)$, where $\bar \lambda \ge 1$ and $\bar x$ satisfies $W \bar x \le w$.
The objective value of~$(\bar \lambda \bar x, \bar \lambda)$ is $\bar \lambda^2 (\bar x_1 - \bar x_2^2 - 2^{-\phi})$.
Since $\bar x$ is feasible to \eqref{pr ParVav}, we have $\bar x_1 - \bar x_2^2 > 2^{-\phi}$, thus the  objective value of~$(\bar \lambda \bar x, \bar \lambda)$ is positive.
In particular, \eqref{pr hard} is bounded.
\end{prf}

In particular, Theorem~\ref{th hard} strengthens the result by Murty and Kabadi~\cite{MurKab87} that deciding whether a QP is bounded or not is NP-hard.

\subsection{Literature review}

In this section, we review the known exact and approximation algorithms for MIQP with a polynomial running time.

MIQP admits a polynomial time approximation algorithm if the dimension $n$ is fixed~\cite{DeLHemKopWei08}.
MIQP is polynomially solvable if the dimension $n$ is fixed and the objective is convex~\cite{Kha83} or concave~\cite{CooHarKanMcDia92,CooKanSch90,HilOerWei15}.
If the objective is concave with a fixed number of negative eigenvalues and the number $p$ of integer variables is fixed, there is a polynomial time approximation algorithm~\cite{dP18}.



Next, we survey \emph{Integer Quadratic Programming} (IQP), which is the special case of MIQP where all variables are integer, \ie $p=n$.
IQP is solvable in polynomial time in dimension one and two~\cite{dPWei14}.
Furthermore, there is a polynomial time approximation algorithm 
if the dimension is fixed and the objective is homogeneous with at most one positive or negative eigenvalue~\cite{HilWeiZem16}.
If the objective function is separable and convex,
and the constraint matrix $W$ is TU, then IQP can be solved in polynomial time~\cite{HocSha90}.
IQP admits a polynomial time approximation algorithm if the objective is separable and concave, with a fixed number of negative eigenvalues, and the largest absolute value of the subdeterminants of the constraint matrix is bounded by two~\cite{dP19}.
Other IQP tractability results under specific structural restrictions can be found in~\cite{LeeOnnRomWei12,EibGanKnoOrd19}.

Finally, we discuss \emph{Quadratic Programming} (QP), the special case of MIQP where all variables are continuous, \ie $p=0$.
QP can be solved in polynomial time if the dimension is fixed~\cite{Vav90,dPDeyMol17}.
Furthermore, QP admits a polynomial time approximation algorithm if the number of negative eigenvalues of~$H$ is fixed~\cite{Vav92i}, and it admits a weaker polynomial time approximation algorithm in general~\cite{Vav93}.
If the objective is convex, then QP can be solved in polynomial time~\cite{KozTarKha79}.


\subsection{Overview of the results and organization of the paper}


Our approximation algorithm is based on the novel concepts of \emph{spherical form MIQP} and of \emph{aligned vectors}.
These two notions significantly enhance the available mathematical toolkit for the design and analysis of algorithms for MIQP, and therefore their importance is not limited to this work.

Section~\ref{sec sym dec} and Section~\ref{sec sform} are devoted to 
finding a change of basis that transforms a MIQP in spherical form.
In a 
\emph{spherical form MIQP} the objective is separable and the polyhedron 
has a ``spherical appearance''.
Moreover, the set $\Z^p \times \R^{n-p}$ is replaced by a set 
of the form $\Lambda + \spn(\Lambda)^\perp$, for a lattice $\Lambda$ of rank $p$. 
The formal definition is given in Section~\ref{sec sform}.
In order to obtain this change of basis we develop a number of results of independent interest.

Since a spherical form MIQP has a separable objective function, in particular we need to find an invertible matrix $L$ and a diagonal matrix $D$ such that $H = LDL^\transp$.
In Section~\ref{sec sym dec} we focus on this simpler task and,
in Theorem~\ref{th decomposition} and Corollary~\ref{cor decomposition}, we present a \emph{symmetric decomposition algorithm} that constructs such matrices $L,D$ in strongly polynomial time.
This is the first known polynomial time algorithm for this problem.

In Section~\ref{sec sform}, we build on this algorithm and obtain, 
in Proposition~\ref{prop simultaneous diagonalization}, 
a rational version of theorems on simultaneous diagonalization of matrices.
In particular, we show that we can find in polynomial time an invertible matrix $L$ that at the same time diagonalizes a given matrix $H$, 
and provides the shape of an ellipsoid that approximates a given polytope within a factor depending only on the dimension.
This result is the main building block that allows us to obtain, in Proposition~\ref{prop change of basis},
a polynomial time algorithm to transform a MIQP in spherical form.


In Section~\ref{sec aligned} we introduce the concept of \emph{aligned vectors} for a spherical form MIQP.
In particular, they are two feasible vectors that are ``far'' in the direction where the objective is ``most curved'' and ``almost aligned'' in all other directions. Furthermore, their midpoint is feasible as well.
We then show, in Proposition~\ref{prop latest}, that if a spherical form MIQP has two aligned vectors, then it is possible to find an $\epsilon$-approximate solution by solving a number of MILPs.
This number is polynomial in~$1/\epsilon$ if both $k$ and $p$ are fixed in the original \eqref{prob main}.

In Section~\ref{sec flatness} 
we focus on the problem of deciding whether a spherical form MIQP has two aligned vectors or not.
In Proposition~\ref{prop aligned} we 
%
give a polynomial time algorithm that either finds two aligned vectors, or finds a vector $v \in \spn(\Lambda)$ along which the polyhedron is ``flat''.
The vector $v$ allows us to decompose the problem in a number of MIQPs with fewer integer variables.
Furthermore, this number depends only on $k$ and $p$, and thus is a constant if both $k$ and $p$ are fixed.

In Section~\ref{sec main algorithm} we then present our approximation algorithm for MIQP and provide a proof of Theorem~\ref{th main}.
The algorithm first uses Proposition~\ref{prop change of basis} to find a change of basis that transforms the input MIQP in spherical form.
Then, it employs Proposition~\ref{prop aligned} and it 
either finds two aligned vectors, or finds a vector $v \in \spn(\Lambda)$ along which the polyhedron is ``flat''.
In the first case, we use Proposition~\ref{prop latest} to find an $\epsilon$-approximate solution.
In the second case, the input MIQP is decomposed into a constant number of instances with fewer integer variable, and the algorithm is recursively applied to these instances.
At the end of the execution, the algorithm returns the best solution found, and we show that it is an $\epsilon$-approximate solution to the input MIQP.
In this paper, we will be using several concepts from computational complexity.
Recall that 
a \emph{strongly polynomial algorithm} is a polynomial space algorithm in the Turing model and a polynomial time algorithm in the arithmetic model.
The definition of strong polynomiality mixes the Turing model and the arithmetic model of computation.
Throughout the paper, unless we state a different model, we mean the Turing model.
For more details on time and space complexity we refer the reader to~\cite{GroLovSch88}.
In particular, we recall that a strongly polynomial algorithm is also a polynomial time algorithm.


\section{A strongly polynomial algorithm for symmetric decomposition}
\label{sec sym dec}

Given a rational symmetric $n \times n$ matrix $\H$, a \emph{symmetric decomposition} of~$\H$ is a decomposition of the form $B\H B^\transp = D$, where $B$ is an $n \times n$ nonsingular matrix and $D$ is an $n \times n$ diagonal matrix.
The goal of this section is to give an algorithm that constructs a symmetric decomposition of any rational symmetric matrix $\H$ with two fundamental properties: 
(i) the algorithm is strongly polynomial,
(ii) the Frobenius norms of~$B$ and $B^{-1}$ are upper bounded by an integer of size polynomial in~$n$.
To the best of our knowledge, our algorithm is the first known polynomial time algorithm that finds a symmetric decomposition of any rational symmetric matrix.
Note that the spectral decomposition, the Schur decomposition, and Takagi's factorization yield a symmetric decomposition of a rational symmetric matrix.
However, none of these decompositions can be performed in polynomial space since the resulting matrices generally contain irrational elements.
Other related matrix decompositions are the Cholesky decomposition and the $LDL^\transp$ decomposition, but are not applicable to indefinite matrices.
For more details on matrix decompositions we refer the reader to~\cite{GolVan4th}.

By introducing pivoting operations that perform symmetric additions of rows and columns, as well as symmetric interchanges, Dax and Kaniel~\cite{DaxKan77} describe an algorithm that constructs a symmetric decomposition of any symmetric $n \times n$ matrix $\H$.
Their algorithm performs a number of arithmetic operations that is polynomial in~$n$, thus it is a polynomial time algorithm in the arithmetic model.
However, it is unknown if it is a polynomial time algorithm or a polynomial space algorithm. 

In this section, we present a strongly polynomial version of Dax and Kaniel's algorithm. 
Therefore, for our version of the algorithm, we show that all numbers stored during the execution of the algorithm have size that is polynomial in the size of the input matrix $\H$.
This in particular implies that the output matrices $B$ and $D$ have polynomial size.
The proof builds on the technique introduced by Edmonds to perform Gaussian elimination in strongly polynomial time~\cite{Edm67}, but it is more involved due to the ``complete pivoting'' performed at each iteration.
In particular, while in Gaussian elimination every number stored during the algorithm is a ratio of subdeterminants of the original matrix, every number stored in our version Dax and Kaniel's algorithm at iteration $k$ is shown to be a ratio of subdeterminants of the matrix obtained from $\H$ by performing only the first $k$ pivoting operations.

Another fundamental property of our symmetric decomposition algorithm is that
the Frobenius norms of~$B$ and $B^{-1}$
are upper bounded by an integer of size polynomial in~$n$.
In particular, this integer depends only on $n$ and not on the input matrix.
This property will be fundamental in the next sections of the paper, where the symmetric decomposition algorithm will be used to obtain a change of basis for our MIQP.
Recall that the \emph{Frobenius norm} of an $m \times n$ matrix $A$ is defined by
$\Fnorm{A} := \sqrt{\sum_{i=1}^m \sum_{j=1}^n A_{ij}^2}.$

Therefore, the purpose of this section is to prove the following result.

\begin{theorem}
\label{th decomposition}
Let $\H$ be a rational symmetric $n \times n$ matrix.
There is a strongly polynomial algorithm that finds matrices $B,D$ such that $B \H B^\transp = D$ is a symmetric decomposition of~$\H$.
Furthermore, 
$\Fnorm{B}$ and $\Fnorm{B^{-1}}$ are upper bounded by $(5n)^{n/2}$.
\end{theorem}

If we set $L:=B^{-1}$ in Theorem~\ref{th decomposition}, we obtain $\H = LDL^\transp$.
Since the inverse can be computed in strongly polynomial time~\cite{Edm67}, also this decomposition can be obtained in strongly polynmomial time.

\begin{corollary}
\label{cor decomposition}
Let $\H$ be a rational symmetric $n \times n$ matrix.
There is a strongly polynomial algorithm that finds
an invertible $n \times n$ matrix $L$ and an $n \times n$ diagonal matrix $D$ such that $\H = L D L^\transp$. 
Furthermore, 
$\Fnorm{L}$ and $\Fnorm{L^{-1}}$ are upper bounded by $(5n)^{n/2}$.
\end{corollary}

Corollary~\ref{cor decomposition} then provides a strongly polynomial algorithm to compute a change of basis that transforms a general \eqref{prob main} in a separable form. 
Namely, compute the decomposition $H = LDL^\transp$ and set $y := L^\transp x$.
In particular, our approach can be substituted to the techniques used in~\cite{Vav92c,HilWeiZem16,dP16,dP18} to transform the original QP or MIQP in a separable form.

In the remainder of this section we only consider matrices that are $n \times n$, thus we avoid repeating it throughout the section.

\subsection{Description of the symmetric decomposition algorithm}
\label{sec sym dec description}

In this section, we describe the symmetric decomposition algorithm that we analyze.
It is the version of Dax and Kaniel's algorithm where the parameter $\gamma$ is always chosen in~$\pm 1$.

Let $\H$ be the rational symmetric $n \times n$ matrix given in the input.
Let $\H^{(0)} := \H$, and for every $k=1,\dots,n-1$, we denote by $\H^{(k)}$ the $n \times n$ matrix obtained after $k$ iterations of the algorithm.
The matrix $\H^{(k)}$ is symmetric and all the off-diagonal elements in the first $k$ rows and columns equal zero.
In particular, $\H^{(n-1)}$ is a diagonal matrix and coincides with the matrix $D$ in the output.

For any $k = 1,\dots,n-1$, we now describe the iteration $k$ of the symmetric decomposition algorithm, where the matrix $\H^{(k)}$ is obtained from $\H^{(k-1)}$.
The iteration is subdivided into two stages, called ``pivoting'' and ``elimination''.

\smallskip
\noindent
\textbf{Pivoting.} 
The goal of the pivoting stage is to ensure that the \emph{pivotal element,} which is the element in the $(k,k)$ position,
is one with largest absolute value among rows and columns $k,\dots,n$.
Let $s$ and $r$ be indices such that 
$\abs{\H^{(k-1)}_{sr}} = \max_{i,j \in \{k,\dots,n\}} \abs{\H^{(k-1)}_{ij}}.$
Since $\H^{(k-1)}_{sr}=\H^{(k-1)}_{rs}$ we can assume without loss of generality that $s \le r$.
Let $\tilde \H$ be the symmetric $n \times n$ matrix obtained from $\H^{(k-1)}$ by interchanging rows $s$ and $k$, and interchanging columns $s$ and $k$.
If $s=r$, then $\H^{(k-1)}_{sr} = \tilde \H_{kk}$.
In this case, we have achieved our goal and the pivoting is terminated.
Thus, we now assume $s < r$.
This implies that $\H^{(k-1)}_{sr} = \tilde \H_{rk}$.
We define
\begin{align*}
\gamma := 
\begin{cases}
+1 & \text{if } \tilde \H_{rk} (\tilde \H_{kk} + \tilde \H_{rr}) \ge 0 \\
-1 & \text{if } \tilde \H_{rk} (\tilde \H_{kk} + \tilde \H_{rr}) < 0,
\end{cases}
\end{align*}
and we let $\doubletilde \H$ be the symmetric $n \times n$ matrix obtained from $\tilde \H$ by adding 
row $r$ multiplied by $\gamma$ to row $k$, and adding column $r$ multiplied by $\gamma$ to column $k$.
It is simple to check that the new $(k, k)$ element is the one with largest absolute value among rows and columns $k,\dots,n$, \ie $\abs{\doubletilde \H_{kk}} = \max_{i,j \in \{k,\dots,n\}} \abs{\doubletilde \H_{ij}}$.

\begin{knownproof}
\begin{observation}
\label{obs max}
We have $\abs{\doubletilde \H_{kk}} = \max_{i,j \in \{k,\dots,n\}} \abs{\doubletilde \H_{ij}}$.
\end{observation}

\begin{prf}
First, we write the elements $\doubletilde \H_{ij}$, for $i,j \in \{k,\dots,n\}$ as functions of the elements of~$\tilde \H$.
\begin{align*}
\doubletilde \H_{ij} & = \tilde \H_{ij} 
&& i,j \neq k \\
\doubletilde \H_{kj} & = \tilde \H_{kj} + \gamma \tilde \H_{rj}
&& 
j \neq k \\
\doubletilde \H_{kk} & = (\tilde \H_{kk} + \gamma \tilde \H_{rk}) + \gamma (\tilde \H_{kr} + \gamma \tilde \H_{rr}) 
= \tilde \H_{kk} + \tilde \H_{rr} + 2 \gamma \tilde \H_{rk},
\end{align*}
where in the last equality 
we used the fact that $\gamma^2 = 1$.
We obtain
\begin{align*}
\abs{\doubletilde \H_{ij}} & \le \abs{\tilde \H_{rk}}
&& i,j \neq k \\
\abs{\doubletilde \H_{kj}} & \le \abs{\tilde \H_{kj}} + \abs{\tilde \H_{rj}} \le 2 \abs{\tilde \H_{rk}}
&& 
j \neq k \\
\abs{\doubletilde \H_{kk}} & 
= \abs{\tilde \H_{kk} + 
\tilde \H_{rr}} + \abs{2 \gamma \tilde \H_{rk}} \ge 2 \abs{\tilde \H_{rk}}.
\end{align*}
In the last equality we used that $\tilde \H_{kk} + \tilde \H_{rr}$ and $2 \gamma \tilde \H_{rk}$ have the same sign due to the definition of~$\gamma$.
This concludes the proof of the observation.
\end{prf}
\end{knownproof}

Pivoting can be achieved via matrix multiplication.
We define the matrix $\tilde P_k$ which interchanges rows $s$ and $k$, thus it is the permutation matrix obtained from the identity matrix by interchanging rows $s$ and $k$
(note that, if $s=k$, then $\tilde P_k$ is the identity matrix).
The matrix $\doubletilde P_k$ adds (if necessary) the row $r$ multiplied by $\gamma$ to row $k$, therefore, it is the identity matrix if $s=r$, or it is obtained from the identity matrix by replacing the zero element in the $(k, r)$ position with the scalar $\gamma$.
The matrix $\tilde \H$ can then be written as
$\tilde \H = \tilde P_k \H^{(k-1)} \tilde P_k^\transp,$
while the matrix $\doubletilde \H$ is the product
$\doubletilde \H  = \doubletilde P_k \tilde \H \doubletilde P_k^\transp = P_k \H^{(k-1)} P_k^\transp,$ where $P_k := \doubletilde P_k \tilde P_k.$

\smallskip
\noindent
\textbf{Elimination.} 
The goal of this stage is to obtain zeros in the off-diagonal elements of row and column $k$.
We first perform row elimination and then column elimination. 
The row elimination is done as in Gaussian elimination:
For each $i=k+1,\dots,n$, add row $k$ multiplied by 
$-\doubletilde \H_{ik} / \doubletilde \H_{kk}$ 
to row $i$. 
The column elimination is done symmetrically:
for each $j=k+1, \dots, n$, add column $k$ multiplied by 
$-\doubletilde \H_{kj} / \doubletilde \H_{kk}$ 
to column $j$. 

Row elimination is performed by multiplying on the left by the matrix $(\id-E_k)$, where $\id$ denotes the $n \times n$ identity matrix and the elements of~$E_k$ are given by 
\begin{align}
\label{eq E elements}
(E_k)_{ik} := 
\doubletilde \H_{ik} / \doubletilde \H_{kk}
\qquad i=k+1,\dots,n, 
\end{align}
and all the other elements are zeros.
Symmetrically, column elimination is performed by multiplying on the right by the matrix $(\id-E_k)^\transp$.
Therefore, the matrix $\H^{(k)}$ is obtained from $\H^{(k-1)}$ via the matrix product
\begin{align}
\label{eq DaxKan 2.1}
\H^{(k)} := [(\id-E_k) P_k] \H^{(k-1)} [(\id-E_k) P_k]^\transp.
\end{align}

\smallskip

This completes the description of the iteration $k$ of the symmetric decomposition algorithm.
At the end of iteration $n-1$ the algorithm returns the diagonal matrix $D := \H^{(n-1)}$ and the nonsingular matrix 
\begin{align*}
B := (\id - E_{n-1}) P_{n-1} \cdots (\id - E_1) P_1.
\end{align*}

It follows directly from the description of the algorithm that the algorithm is correct, \ie $B\H B^\transp = D$ is a symmetric decomposition of~$\H$.

\subsection{Analysis of the algorithm}

In this section, we prove the first part of Theorem~\ref{th decomposition}.
Namely, we show that the symmetric decomposition algorithm presented in Section~\ref{sec sym dec description} runs in strongly polynomial time.
Clearly, the number of arithmetic operations performed is polynomial in~$n$.
Therefore, we only need to show that the size of each matrix constructed during the execution is polynomial in the size of~$\H$.
For matrices $\tilde P_k, \doubletilde P_k, P_k$, for $k = 1,\dots,n-1$, this follows directly from their definition.
In fact, we only need to show that each matrix $\H^{(k)}$, for $k = 1,\dots,n-1$, has size polynomial in the size of~$\H$.
Indeed, once this is proven, we obtain that also $E_k$ and the returned matrix $B$ have size polynomial in the size of~$\H$.

Thus, we now focus our attention on the matrix $\H^{(k)}$.
From the equality \eqref{eq DaxKan 2.1} we deduce that
\begin{align*}
\H^{(k)} = B^{(k)} \H {B^{(k)}}^\transp,
\end{align*}
where
\begin{align*}
B^{(k)} := 
(\id - E_k) P_k
(\id - E_{k-1}) P_{k-1} 
\cdots (\id - E_1) P_1.
\end{align*}
As noticed on page 224 in~\cite{DaxKan77}, 
it is simple to verify that for every $t,j \in \{1,\dots, n-1\}$ with $t < j$, we have $E_t P_j = E_t$.
This in turn implies that for every \mbox{$t,j \in \{1,\dots, n-1\}$} with $t < j$, we have
\begin{align*}
P_j (\id - P_{j-1} P_{j-2} \cdots P_{t+1} E_t) = (\id - P_j P_{j-1} \cdots P_{t+1} E_t) P_j,
\end{align*}
which allows us to write $B^{(k)}$ in the form
\begin{align}
\label{eq Bk decomposition}
B^{(k)} & = (\id - E_k) (\id - P_k E_{k-1}) 
\cdots (\id - P_k \cdots P_2 E_1) P_k \cdots  P_1.
\end{align}
Therefore $\H^{(k)}$ can be written as 
\begin{align}
\label{eq new eli}
\H^{(k)} & = E^{(k)} \G^{(k)} {E^{(k)}}^\transp,
\end{align}
where 
\begin{align*}
\G^{(k)} & := (P_k \cdots P_1) \H (P_k \cdots P_1)^\transp, \\
E^{(k)} & := (\id - E_k) (\id - P_k E_{k-1}) 
\cdots (\id - P_k \cdots P_2 E_1).
\end{align*}

In the next lemma, we analyze the matrices $P_{k} P_{k-1} \cdots P_{t+1} E_t$ in the definition of~$E^{(k)}$. 
The second part of the statement will only be used later in Section~\ref{sec sym dec Frobenius}.

\begin{lemma}
\label{lem PsE}
For each $t \in \{1,\dots,k\}$, the matrix $P_k P_{k-1} \cdots P_{t+1} E_{t}$ can have nonzeros only in positions $(t+1,t),\dots,(n,t)$.
Furthermore, the elements in rows $t+1,\dots,k$ are bounded by two in absolute value, while the elements in rows $k+1,\dots,n$ are bounded by one in absolute value.
\end{lemma}

\begin{prf}
We show this lemma by induction on $k-t$.
In the base case we have $k-t = 0$, thus we are considering matrix $E_t$.
By definition, $E_t$ can have nonzeros only in positions $(t+1,t),\dots,(n,t)$, and from \eqref{eq E elements} all nonzeros are are bounded by one in absolute value.

For the inductive step we assume $k-t \ge 1$ and consider the matrix $P_k ( P_{k-1} \cdots P_{t+1} E_{t})$.
By induction, $P_{k-1} \cdots P_{t+1} E_{t}$ can have nonzeros only in positions $(t+1,t),\dots,(n,t)$.
Furthermore, the elements in rows $t+1,\dots,k-1$ are bounded by two in absolute value, while the elements in rows $k,\dots,n$ are bounded by one in absolute value.
We have $P_k = \doubletilde P_k \tilde P_k,$ where the matrix $\tilde P_k$ interchanges two rows in~$\{k,\dots,n\}$,
and the matrix $\doubletilde P_k$ adds or subtracts (if necessary) a row in~$\{k+1,\dots,n\}$
to row $k$. 
Since $k \ge t+1$, the matrix $P_k ( P_{k-1} \cdots P_{t+1} E_{t})$ can have nonzeros only in positions $(t+1,t),\dots,(n,t)$.
The elements in rows $t+1,\dots,k-1$ are left unchanged, thus they are bounded by two in absolute value.
The element in row $k$ is now bounded by two in absolute value, while the elements in rows $k+1,\dots,n$ remain bounded by one in absolute value.
\end{prf}

Next, we use Lemma~\ref{lem PsE} to discuss the effect of multiplying a matrix on the left by $E^{(k)}$.
Note that a multiplication of this type is performed in \eqref{eq new eli}.

\begin{lemma}
\label{lem Ek left}
Multiplying a matrix on the left by $E^{(k)}$ results in a sequence of elementary row operations in which a multiple of a row $t \in \{1,\dots,k\}$ is added to a row in~$\{t+1, \dots, n\}$.
\end{lemma}

\begin{prf}
Due to the definition of~$E^{(k)}$, it suffices to show that multiplying a matrix on the left by
$(\id - P_k \cdots P_{t+1} E_{t}),$
for $t \in \{1,\dots,k\}$,
results in a sequence of elementary row operations in which a multiple of row $t$ is added to a row in~$\{t+1,\dots, n\}$.

From Lemma~\ref{lem PsE}, the matrix $P_k \cdots P_{t+1} E_{t}$ can have nonzeros only in positions $(t+1,t),\dots,(n,t)$.
Hence, the multiplication on the left by matrix $(\id - P_k \cdots P_{t+1} E_{t})$ preserves the first $t$ rows, and each subsequent row is obtained by adding a multiple of row $t$ to the original row.
\end{prf}

We are finally ready to show, in the next claim, that each matrix $\H^{(k)}$ has size polynomial in the size of~$\H$.
This concludes the proof that our symmetric decomposition algorithm runs in strongly polynomial time.

\begin{claim}
For each $k \in \{1,\dots,n-1\}$, the size of~$\H^{(k)}$ is polynomial in the size of~$\H$.
\end{claim}

\begin{prf}
Let $\G^{(k)}_k$ denote the $k \times k$ submatrix of~$\G^{(k)}$ determined by the first $k$ rows and columns, and let $\G^{(k)}_{k\mid ij}$, for $i,j \in \{k+1,\dots,n\}$, denote the $(k+1) \times (k+1)$ submatrix of~$\G^{(k)}$ determined by rows $\{1,\dots,k,i\}$ and columns $\{1,\dots,k,j\}$.

It suffices to show that for every $k \in \{1,\dots,n-1\}$ and for every $i,j \in \{k+1,\dots,n\}$, we have
\begin{align}
\label{eq Gauss}
\H_{ij}^{(k)} = 
\det(\G^{(k)}_{k\mid ij}) / \det(\G^{(k)}_k).
\end{align}
In fact, the definition of~$\G^{(k)}$ implies that its size is polynomial in the size of~$\H$.
Therefore, also $\det(\G^{(k)}_{k\mid ij})$ and $\det(\G^{(k)}_k)$ have size polynomial in the size of~$\H$, and so does each element of~$\H^{(k)}$ due to \eqref{eq Gauss}.
Therefore, in the remainder of the proof we show \eqref{eq Gauss}.

Consider the product $E^{(k)} \G^{(k)}$ in \eqref{eq new eli}.
From Lemma~\ref{lem Ek left}, the resulting matrix is obtained from $\G^{(k)}$ via a sequence of elementary row operations.
Among all these elementary row operations, only a subset modify the first $k$ rows of the matrix $\G^{(k)}$.
From Lemma~\ref{lem Ek left}, in each of these elementary row operations, a multiple of a row $t \in \{1,\dots,k-1\}$ is added to a row in~$\{t+1,\dots,k\}$.
Similarly, among all the elementary column operations performed by ${E^{(k)}}^\transp$, only a subset modify the first $k$ columns of the matrix $\G^{(k)}$.
In each of these elementary column operations, a multiple of a column $t \in \{1,\dots,k-1\}$ is added to a column in~$\{t+1,\dots,k\}$.
We perform these subsets of elementary row and column operations to the matrix $\G^{(k)}_k$.
From \eqref{eq new eli}, the resulting matrix is precisely the submatrix of~$\H^{(k)}$ given by the first $k$ rows and columns, hence it is diagonal with elements $\H_{11}^{(k)}, \dots ,\H_{kk}^{(k)}$ in the diagonal.
Note that each elementary operation considered preserves the determinant of~$\G^{(k)}_k$.
Thus,
\begin{align}
\label{eq den}
\det(\G^{(k)}_k) = \H_{11}^{(k)} \cdots \H_{kk}^{(k)}.
\end{align}

A similar argument can be applied to the matrix $\G^{(k)}_{k\mid ij}$.
Among all the elementary row operations performed by $E^{(k)}$, only a subset modify rows $\{1,\dots,k,i\}$ of the matrix $\G^{(k)}$. 
From Lemma~\ref{lem Ek left}, in each of these elementary row operations, a multiple of a row $t \in \{1,\dots,k\}$ is added to a row in~$\{t+1,\dots,k,i\}$.
Similarly, among all the elementary column operations performed by ${E^{(k)}}^\transp$, only a subset modify columns $\{1,\dots,k,j\}$ of the matrix $\G^{(k)}$.
In each of these elementary column operations, a multiple of a column $t \in \{1,\dots,k\}$ is added to a column in~$\{t+1,\dots,k,j\}$.
We perform this subset of elementary operations to the matrix $\G^{(k)}_{k\mid ij}$.
From \eqref{eq new eli}, the resulting matrix is precisely the submatrix of~$\H^{(k)}$ determined by rows $\{1,\dots,k,i\}$ and columns $\{1,\dots,k,j\}$.
Hence it is diagonal with elements $\H_{11}^{(k)}, \dots ,\H_{kk}^{(k)},\H_{ij}^{(k)}$ in the diagonal.
Each elementary operation considered preserves the determinant of~$\G^{(k)}_{k\mid ij}$.
Thus, we have
\begin{align*}
\det(\G^{(k)}_{k\mid ij}) = \H_{11}^{(k)} \cdots \H_{kk}^{(k)} \cdot \H_{ij}^{(k)}.
\end{align*}
Dividing the latter equation by equation \eqref{eq den}, we obtain \eqref{eq Gauss}.
\end{prf}

\subsection{Frobenius norm of~$B$ and $B^{-1}$}
\label{sec sym dec Frobenius}

In this section, we prove the second part of Theorem~\ref{th decomposition}.
Namely, we show that for the matrix $B$ returned by the algorithm described in Section~\ref{sec sym dec description}, both its Frobenius norm and the Frobenius norm of its inverse are upper bounded by an integer of size polynomial in~$n$.
We will use the fact that the Frobenius norm is \emph{submultiplicative}, \ie for matrices $A,A'$ we have $\Fnorm{AA'} \le \Fnorm{A} \cdot \Fnorm{A'}$.

\begin{claim}
\label{claim Fr norms}
The Frobenius norm of matrices $B$ and $B^{-1}$ is upper bounded by $(5n)^{n/2}$.
\end{claim}

\begin{prf}
From \eqref{eq Bk decomposition},
we can write $B = B^{(n-1)}$ as the product $B = EP$, where 
\begin{align*}
E & := E^{(n-1)} = (\id - E_{n-1}) (\id - P_{n-1} E_{n-2}) 
\cdots (\id - P_{n-1} \cdots P_2 E_1), \\
P & := P_{n-1} P_{n-2} \cdots P_{1}.
\end{align*}
In order to bound the Frobenius norm of~$B$ and $B^{-1}$, we bound separately the Frobenius norm of~$P, P^{-1}, E$, and $E^{-1}$.

\smallskip
\noindent
\ul{Norm of~$P$.}
Recall that each matrix $P_k$ is the product of two matrices $P_k = \doubletilde{P}_k \tilde{P}_k$, where
the matrix $\tilde P_k$ interchanges row $s$ and $k$, where $s \ge k$, and the matrix $\doubletilde P_k$ adds (if necessary) row $r$ multiplied by $\gamma$ to row $k$, where $r > k$.
Therefore, for each $k=1,\dots,n-1$, in the matrix $P_k \cdots P_1$, the last $n-k$ rows are permuted rows of the identity matrix, while each of the first $k$ rows has at most two nonzero elements, each one being $\pm 1$.
We obtain 
\begin{align*}
\Fnorm{P} & = \Fnorm{P_{n-1} \cdots P_1} \le \sqrt{2n-1}.
\end{align*}

\smallskip
\noindent
\ul{Norm of~$P^{-1}$.}
Each matrix $P_k^{-1}$ is the product $P_k^{-1} = \tilde{P}_k^{-1} \doubletilde{P}_k^{-1}$, where the matrix $\doubletilde P_k^{-1}$ adds (if necessary) row $r$ multiplied by $-\gamma$ to row $k$, where $r > k$, and the matrix $\tilde P_k^{-1}$ interchanges row $s$ and $k$, where $s \ge k$.
Therefore, for each \mbox{$k=n-1,\dots,1$}, in the matrix $P_k^{-1} P_{k+1}^{-1} \cdots P_{n-1}^{-1}$, the first $k-1$ rows coincide with the first $k-1$ rows of the identity matrix, and the remaining rows are a permutation of the rows from $k$ to $n$ of an upper triangular matrix with elements in~$0, \pm 1$.
We obtain
\begin{align*}
\Fnorm{P^{-1}} = \Fnorm{P_1^{-1} P_2^{-1} \cdots P_{n-1}^{-1}} 
\le 
\sqrt{(n^2 + n)/2}.
\end{align*}

\smallskip
\noindent
\ul{Norm of~$E$.}
From Lemma~\ref{lem PsE} with $k=n-1$, for each $t \in \{1,\dots,n-1\}$, the matrix $P_{n-1} P_{n-2} \cdots P_{t+1} E_t$ can have nonzeros only in positions $(t+1,t),\dots,(n,t)$.
Furthermore, the elements in rows $t+1,\dots,n-1$ are bounded by two in absolute value, while the element in the last row is bounded by one in absolute value.
Thus, we obtain 
\begin{align*}
\Fnorm{\id-P_{n-1} P_{n-2} \cdots P_{t+1} E_t} \le \sqrt{(n+1) + 4(n-t-1)} \le \sqrt{5(n-1)}.
\end{align*}
Hence
\begin{align*}
\Fnorm{E} & 
\le \Fnorm{\id - E_{n-1}} \cdot \Fnorm{\id - P_{n-1} E_{n-2}} \cdots \Fnorm{\id - P_{n-1} \cdots P_2 E_1} \\
& \le \sqrt{(5(n-1))^{n-1}}.
\end{align*}

\smallskip
\noindent
\ul{Norm of~$E^{-1}$.}
Once again, from Lemma~\ref{lem PsE} with $k=n-1$, we know that for each $t \in \{1,\dots,n-1\}$, the matrix $P_{n-1} P_{n-2} \cdots P_{t+1} E_t$ can have nonzeros only in positions $(t+1,t),\dots,(n,t)$.
This fact allows us to write $E^{-1}$ as
\begin{align*}
E^{-1} & = (\id - P_{n-1} \cdots P_2 E_1)^{-1} \cdots (\id - P_{n-1} E_{n-2})^{-1} (\id - E_{n-1})^{-1} \\
& = (\id + P_{n-1} \cdots P_2 E_1) \cdots (\id + P_{n-1} E_{n-2}) (\id + E_{n-1}) \\
& = \id + P_{n-1} \cdots P_2 E_1 + \cdots + P_{n-1} E_{n-2} + E_{n-1}.
\end{align*}
In particular, the matrix $E^{-1}$ is \emph{unit lower triangular,} \ie lower triangular with all elements on the main diagonal equal to one.
The second part of Lemma~\ref{lem PsE} then implies that the elements in rows $1,\dots,n-1$ are bounded by two in absolute value, while the elements in the last row are bounded by one in absolute value.
We obtain 
\begin{align*}
\Fnorm{E^{-1}} \le \sqrt{(2n-1) + 4 (n^2 - 3n +2)/2} = \sqrt{2 n^2 - 4n + 3}.
\end{align*}

\smallskip
\noindent
\ul{Norm of~$B$ and $B^{-1}$.}
Using the obtained bounds on the Frobenius norm of~$P$, $P^{-1}$, $E$, $E^{-1}$, we derive
\begin{align*}
\Fnorm{B} & = \Fnorm{EP} \le \Fnorm{E} \, \Fnorm{P} 
\le \sqrt{(5(n-1))^{n-1} (2n-1)}, \\
\Fnorm{B^{-1}} & = \Fnorm{P^{-1} E^{-1}} \le \Fnorm{P^{-1}} \, \Fnorm{E^{-1}} 
\le 
\sqrt{(n^2 + n) (2 n^2 - 4n + 3)/2}.
\end{align*}
It can be checked that $(5n)^{n/2}$ is larger than both upper bounds for any $n$.
Therefore, the claim follows.
\end{prf}

While the bound on the Frobenius norm of matrices $B$ and $B^{-1}$ in Claim~\ref{claim Fr norms} is sufficient for our task, we remark that a better bound can be obtained by providing a better bound on $\Fnorm{E}$.
This can be done by bounding the largest absolute value of an element in~$E$, instead of using the fact that the Frobenius norm is submultiplicative.


\section{Simultaneous diagonalization and spherical form MIQP}
\label{sec sform}

In this section, a fundamental role is played by the spherical form MIQP.
To formally define this problem, we now briefly recall the notion of lattice, and introduce some notation.

Given linearly independent vectors 
$b^1, \dots , b^p$
in~$\R^\did$, the \emph{lattice} generated by 
$b^1, \dots , b^p$
is the set
$\Lambda 
 := \bra{\sum_{i=1}^p v_i b^i : v_i \in \Z \ \forall i=1,\dots,p}$
of integer linear combinations of the vectors $b^i$, for $i=1,\dots,p$.
The \emph{rank} of~$\Lambda$ is $p$ and the \emph{dimension} of~$\Lambda$ is $\did$. 
If $p = \did$, then $\Lambda$ is said to be a \emph{full rank lattice}.
Note that, in this paper, we will consider mainly lattices that are not full rank.
The vectors 
$b^1, \dots , b^p$
are called a \emph{lattice basis}.
%
Given a vector $\ctr \in \R^\did$ and a nonnegative scalar $r$, we denote by $\B(\ctr,r)$ the \emph{closed ball} with center $\ctr$ and radius $r$. 
Formally, 
\begin{align*}
\B(\ctr,r) := \{x \in \R^\did : \norm{x-\ctr} \le r\}.
\end{align*}
Note that, throughout the paper, we use the \emph{euclidian vector norm} defined
as $\norm{x} := \sqrt{x^\transp x}.$
Given vectors $x^1,\dots,x^t$, we denote by $(x^1,\dots,x^t)$ the vector $({x^1}^\transp,\dots,{x^t}^\transp)^\transp$.
The orthogonal complement of a linear space $\L$ is denoted by $\L^\perp$.

We are now in a position to give the formal definition of a spherical form MIQP.
A \emph{spherical form MIQP} is an optimization problem of the form
\begin{align}
\label{prob s main}
\tag{S-MIQP}
\begin{split}
\min & \quad y^\transp D y + c^\transp y + \zcoeff^\transp z  \\
\st & \quad (y,z) \in \P  \\
& \quad y \in \Lambda + \spn(\Lambda)^\perp, \ z \in \R^{n-\did}.
\end{split}
\end{align}
In this formulation, the variables are $y \in \R^\did$ and $z \in \R^{n-\did}$.
The matrix $D \in \Q^{\did \times \did}$ is diagonal and its diagonal elements satisfy $\abs{D_{11}} \ge \cdots \ge \abs{D_{\did\did}}$. 
Furthermore, $c \in \Q^\did$, and $\zcoeff \in \Q^{n-\did}$.
The set $\Lambda$ is a lattice of rank $p$ and dimension $\did$,
and is given via a rational lattice basis.
Finally, the set $\P \subseteq \R^n$ is a polytope given via a finite system of rational linear inequalities, and it satisfies 
\begin{align}
\label{eq containment}
\B(\ctr,1) \subset \proj_y \P \subset \B(\ctr,\constelli_\did), 
\end{align}
where 
$\proj_y \P$ denotes the orthogonal projection of~$\P$ onto the space $\R^\did$ of the $y$ variables,
$a$ is a given vector in~$\Q^\did$, and
$\constelli_\did$ 
is an integer of size polynomial in~$\did$. 

The symmetric decomposition algorithm described in Section~\ref{sec sym dec} allows us to obtain,  in strongly polynomial time, a change of basis that directly transforms \eqref{prob main} in a separable form.
In this section, our main goal is to obtain another change of basis that not only maps \eqref{prob main} in a separable form, but 
also guarantees that the resulting problem is in spherical form.
The additional requirements on the change of basis will result in an algorithm that is polynomial time instead of strongly polynomial.
To obtain this change of basis, we rely on two key results: (i) the symmetric decomposition algorithm discussed in Section~\ref{sec sym dec}, and 
(ii) the existence of an algorithm based on linear programming that, for every full-dimensional polytope~$\P$, constructs a pair of concentric ellipsoids $\E_1$, $\E_2$ such that $\E_1 \subset \P \subset \E_2$ and $\E_1$ is obtained by shrinking $\E_2$ by a factor depending only on the dimension~\cite{Len83}.

\subsection{Simultaneous diagonalization}

The first result of this section does not deal directly with MIQP but is the main building block that will allow us to transform a MIQP in spherical form.
This result can be interpreted as a rational version of classic theorems on simultaneous diagonalization of matrices (see Section~8.7 in~\cite{GolVan4th}).

In order to present our result we need to introduce ellipsoids.
An \emph{ellipsoid} in~$\R^n$ is an affine transformation of the unit ball. 
That is, an ellipsoid is a set of the form 
\begin{align*}
\E(\ctr,L) = \{x \in \R^n : \norm{L^\transp(x-\ctr)} \le 1\},
\end{align*}
where $\ctr \in \R^n$ and $L$ is an $n \times n$ invertible matrix.
Note that
$
\B(\ctr,r) = \E(\ctr,\id_n/r),
$
where $\id_n$ denotes the $n \times n$ identity matrix.

In what follows, we will often work with rational linear subspaces. 
In the context of polynomial time algorithms, it is not important if they are given to us via a system of linear equations or via a basis, since each description can be obtained in polynomial time from the other.
Given a linear subspace $\L$ of~$\R^n$ of dimension~$\did$, a \emph{basis matrix} of~$\L$ is an $n \times \did$ matrix whose columns $b^1,\dots,b^\did$ form a basis of~$\L$.
An \emph{$\L$-ellipsoid} is a set of the form 
\begin{align*}
\E_\L(\ctr,L) = \{x \in \L : \norm{L^\transp(x-\ctr)} \le 1\},
\end{align*}
where $\L$ is a linear subspace of~$\R^n$,
$\ctr \in \L$,
and $L$ is a basis matrix of~$\L$.
Given a linear subspace $\L$ of~$\R^n$ and a set $\S \subseteq \R^n$, we denote by $\proj_\L (\S)$ the orthogonal projection of~$\S$ onto $\L$.
We also say that a polyhedron $\{x : Wx \le w\}$ is \emph{rational} if $W$ and $w$ are rational.
We are now ready to present the first result of this section.

\begin{proposition}
\label{prop simultaneous diagonalization}
Let $H$ be a rational symmetric $n \times n$ matrix of rank $k$,
let $\{x \in \R^n: Wx \le w\}$ be a full-dimensional rational polytope,
and let $\M$ be a rational linear subspace of~$\R^n$ of dimension $p$.
There is a polynomial time algorithm that finds a linear subspace $\L$ of~$\R^n$ 
containing $\M$ and of dimension $\did$ with $\max\{k,p\} \le \did \le  k+p$, 
a $\did \times \did$ diagonal matrix $D$,
and an $\L$-ellipsoid $\E_{\L} (\ctr,L)$ such that
\begin{itemize}
\item[(i)]
$H = L D L^\transp$,
\item[(ii)]
$\E_{\L} (\ctr, L) \subset \proj_{\L} \{x: Wx \le w\} \subset \E_{\L} (\ctr,L / (2 \did^{3/2} \lceil (5d)^{d/2}\rceil^2)).$ 
\end{itemize}
\end{proposition}

\begin{prfh}
By Corollary~\ref{cor decomposition} there is a strongly polynomial algorithm that computes 
an invertible $n \times n$ matrix $L_1$ and an $n \times n$ diagonal matrix $D_1$ such that $H = L_1 D_1 L_1^\transp$. 
Since $H$ has rank $k$ and $L_1$ is invertible, the matrix $D_1$ has also rank $k$.
Let $D_2$ be the matrix obtained from $D_1$ by deleting row $i$ and column $i$ for each $i$ such that the $i$th diagonal element of~$D_1$ is zero.
Clearly, $D_2$ is an invertible $k \times k$ diagonal matrix.
We also define the matrix $L_2$, obtained from $L_1$ by deleting column $i$ for each $i$ such that the $i$th diagonal element of~$D_1$ is zero.
The matrix $L_2$ is then an $n \times k$ matrix of rank $k$.
Since row and column $i$ of~$D_1$ have all zero elements, we have $H = L_1 D_1 L_1^\transp = L_2 D_2 L_2^\transp$.

Let $\L$ be the linear subspace of~$\R^n$ obtained as the Minkowski sum of~$\M$ and of the linear space spanned by the $k$ columns of~$L_2$.
Clearly, $\L$ contains $\M$, and its dimension $\did$ satisfies $\max\{k,p\} \le \did \le  k+p$.
Note that $\proj_{\L} \{x: Wx \le w\}$ is full-dimensional.
It then follows form Sections~2 and 5 in~\cite{Len83} that there is a polynomial time algorithm which computes an $\L$-ellipsoid $\E_{\L} (\ctr,C)$ such that
\begin{align}
\label{eq elli cont}
\E_{\L}(\ctr, C) \subset \proj_{\L} \{x: Wx \le w\} \subset \E_{\L}(\ctr,C / (2 \did^{3/2})).
\end{align}

Since the $n \times \did$ matrix $C$ is a basis matrix of~$\L$ and each column of~$L_2$ is a vector in~$\L$, we can compute in polynomial time a $\did \times k$ matrix $M$ such that $L_2 = C M$.
We obtain 
\begin{align*}
H = L_2 D_2 L_2^\transp = C M D_2 M^\transp C^\transp = C \tilde H C^\transp,
\end{align*}
where $\tilde H := M D_2 M^\transp$ is a $\did \times \did$ symmetric matrix.

By Corollary~\ref{cor decomposition}, applied to $\tilde H$, there is a strongly polynomial algorithm which computes an invertible $\did \times \did$ matrix $\tilde L$ and a $\did \times \did$ diagonal matrix $\tilde D$ such that $\tilde H = \tilde L \tilde D \tilde L^\transp$. 
Furthermore, $\Fnorm{\tilde L}$ and $\Fnorm{\tilde L^{-1}}$ are upper bounded by $\constdec_\did:= \ceil{(5d)^{d/2}}$.
Note that $\constdec_\did$ is an integer of size polynomial in~$\did$.
We obtain
$H = C \tilde L \tilde D \tilde L^\transp C^\transp.$
By defining the $\did \times \did$ matrix $D$ and the $n \times \did$ matrix $L$ in the statement as
$D := \tilde D / \constdec_\did^2,$ 
$L := \constdec_\did C \tilde L,$
we obtain $H = L D L^\transp$. 
Clearly, $D$ is diagonal, thus condition (i) in the statement holds.

Note that the vector $\ctr$ is in~$\L$.
Moreover, since $C$ is a basis matrix of~$\L$ and $\tilde L$ is invertible, we have that also $L$ is a basis matrix of~$\L$.
Hence $\E_{\L} (\ctr,L)$ is an $\L$-ellipsoid.
We now show that condition (ii) is satisfied.
Using the fact that the Frobenius norm is submultiplicative and that $\Fnorm{\tilde L}$ and $\Fnorm{\tilde L^{-1}}$ are upper bounded by $\constdec_\did$, we obtain
\begin{align*}
\norm{C^\transp (x - \ctr)} &= \norm{\tilde L^{-\transp} L^\transp (x - \ctr)}/\constdec_\did \le \Fnorm{\tilde L^{-1}} \, \norm{L^\transp (x - \ctr)} / \constdec_\did \le \norm{L^\transp (x - \ctr)}, \\
\norm{L^\transp (x - \ctr)} & = \constdec_\did \norm{\tilde L^\transp C^\transp (x - \ctr)} \le \constdec_\did \Fnorm{\tilde L} \, \norm{C^\transp (x - \ctr)} \le \constdec_\did^2 \norm{C^\transp (x - \ctr)}.
\end{align*}
The first chain of inequalities and \eqref{eq elli cont} imply
\begin{align*}
\E_{\L} (\ctr, L) \subseteq \E_{\L} (\ctr, C) \subset \proj_{\L} \{x: Wx \le w\}.
\end{align*}
The second chain of inequalities implies $\E_{\L} (\ctr,\constdec_\did^2 C) \subseteq \E_{\L} (\ctr, L)$, thus from \eqref{eq elli cont},
\begin{equation*}
\proj_{\L} \{x: Wx \le w\} \subset \E_{\L} (\ctr, C / (2 \did^{3/2})) \subseteq \E_{\L} (\ctr, L / (2 \did^{3/2} \constdec_\did^2)). \qedhere
\end{equation*}
\end{prfh}

Consider now the simplest case of Proposition~\ref{prop simultaneous diagonalization}, where we set $\M := \R^n$.
Then $\L = \R^n$, $d=n$, the $\L$-ellipsoids are just ellipsoids, and the polytope $\proj_{\L} \{x: Wx \le w\}$ is simply $\{x: Wx \le w\}$.
In this case, Proposition~\ref{prop simultaneous diagonalization} provides a matrix $L$ that at the same time diagonalizes $H$ and provides the shape of an ellipsoid that approximates the given polytope within a factor depending only on the dimension.
This special case 
can then be interpreted as a rational version of theorems on simultaneous diagonalization of matrices.
If we perform the change of basis $y := L^\transp x$, the given matrix $H$ is diagonalized, and the ellipsoids are just balls.

\subsection{Reduction to spherical form MIQP}

Next, we employ Proposition~\ref{prop simultaneous diagonalization} to show that \eqref{prob main} can be transformed in spherical form \eqref{prob s main}.
Throughout the paper, we denote by $e^1,e^2\dots,e^n$ the standard basis of~$\R^n$.


\begin{proposition}
\label{prop change of basis}
Consider \eqref{prob main}, assume that $\{x: Wx \le w\}$ is a full-dimensional polytope, and let $k$ denote the rank of~$H$.
There is a polynomial time algorithm that finds a change of basis that transforms \eqref{prob main} in spherical form \eqref{prob s main}, where $\did$ satisfies $\max\{k,p\} \le \did \le k+p$,
the rank of the matrix $D$ is $k$,
and 
$\constelli_\did$ in \eqref{eq containment} is the ceiling of~$2 \did^{3/2} \ceil{(5d)^{d/2}}^2$.
\end{proposition}

\begin{prf} 
Consider \eqref{prob main}, assume that $\{x: Wx \le w\}$ is a full-dimensional polytope, and let $k$ denote the rank of~$H$.
By Proposition~\ref{prop simultaneous diagonalization} with $\M := \R^p \times \{0\}^{n-p}$, we obtain in polynomial time a linear subspace $\L$ of~$\R^n$ containing $\M$ and of dimension $\did$ with $\max\{k,p\} \le \did \le  k+p$,
a $\did \times \did$ diagonal matrix $D$, 
and an $\L$-ellipsoid $\E_{\L}(\ctr, L_y)$ such that 
$H = L_y D L_y^\transp$ and
\begin{align}
\label{eq prop}
\E_{\L} (\ctr, L_y) \subset \proj_{\L} \{x: Wx \le w\} \subset \E_{\L} (\ctr, L_y / \constelli_\did),
\end{align}
where we define $\constelli_\did$ as the ceiling of~$2 \did^{3/2} \ceil{(5d)^{d/2}}^2$.
Since $L_y$ is an $n \times \did$ matrix of rank $\did$, it is simple to check that 
the rank of~$D$ coincides with the rank of~$H$.

We now compute an $n \times (n-\did)$ basis matrix $L_z$ of the orthogonal complement $\L^\perp$ of~$\L$.
Denote by $L$ the $n \times n$ invertible matrix $(L_y \mid L_z)$.
We perform the change of basis $x \mapsto (y,z)$, where $(y,z) \in \R^n$ is defined by $(y,z) := L^\transp x$, \ie $y \in \R^\did$ is defined by $y := L_y^\transp x$, and $z \in \R^{n-\did}$ is defined by $z := L_z^\transp x$.

Next, we consider the problem obtained from \eqref{prob main} via the above change of basis, and we show that it coincides with \eqref{prob s main}.
The objective function of the new problem is
\begin{align*}
x^\transp H x + h^\transp x = x^\transp L_y D L_y^\transp x + h^\transp x = y^\transp D y + h^\transp L^{-\transp} (y,z),
\end{align*}
which coincides with the objective function of \eqref{prob s main} if we define the vectors $c \in \Q^\did$ and $\zcoeff \in \Q^{n-\did}$ by 
$(c,\zcoeff) := L^{-1} h$.
The image of the polytope $\{x: Wx \le w\}$ is the set 
$\P := \{(y,z) : WL^{-\transp} (y,z) \le w\}$.
Clearly, $\P$ is a polytope defined by a finite system of rational linear inequalities.

By definition of~$L_z$, the linear subspace $\L$ can be written as $\L = \{x : L_z^\transp x = 0\}$, thus the image of~$\L$ under the change of basis is $\{(y,z) : z = 0\} = \R^\did \times \{0\}^{n-\did}$.
Similarly, the linear subspace $\L^\perp$ can be written as $\L^\perp = \{x : L_y^\transp x = 0\}$, thus the image of~$\L^\perp$ is $\{0\}^\did \times \R^{n-\did}$.

Next we show that \eqref{eq prop} implies \eqref{eq containment}.
The above discussion implies that a point $\proj_{\L} (x)$ is mapped to $\proj_y (L^{-\transp} (y,z)) \times \{0\}^{n-\did}$.
Thus, $\proj_{\L} \{x : Wx \le w\}$ is mapped to $\proj_y \P \times \{0\}^{n-\did}$.
The set $\E_{\L} (\ctr, L_y)$ is mapped to 
\begin{align*}
\bra{(y,z) : \norm{y-L_y^\transp \ctr} \le 1} \cap (\R^\did \times \{0\}^{n-\did}) 
& = \E (L_y^\transp \ctr, \id_\did) \times \{0\}^{n-\did} \\
& = \B(L_y^\transp \ctr,1) \times \{0\}^{n-\did}.
\end{align*}
Similarly, the set $\E_{\L} (\ctr, L_y/\constelli_\did)$ is mapped to 
\begin{align*}
\bra{(y,z) : \norm{(y-L_y^\transp \ctr)/\constelli_\did} \le 1} \cap (\R^\did \times \{0\}^{n-\did}) 
& = \E (L_y^\transp \ctr, \id_\did / \constelli_\did) \times \{0\}^{n-\did} \\
& = \B(L_y^\transp \ctr,\constelli_\did) \times \{0\}^{n-\did}.
\end{align*}
From \eqref{eq prop}, we obtain 
$\B(L_y^\transp \ctr,1) \subset \proj_y \P \subset \B(L_y^\transp \ctr,\constelli_\did),$
which coincides with \eqref{eq containment} if we redefine the vector $\ctr \in \Q^\did$ to be $L_y^\transp \ctr$.

We now consider the image of~$\Z^p \times \R^{n-p}$.
The set $\Z^p \times \R^{n-p}$ can be written as the Minkowski sum $(\Z^p \times \{0\}^{n-p}) + \N + \L^\perp$, where $\N$ is the orthogonal complement of~$\M$ in~$\L$.
Since $\M \subseteq \L$ and the image of~$\L$ is $\R^\did \times \{0\}^{n-\did}$,
we have that the image of~$\Z^p \times \{0\}^{n-p}$ is $\Lambda \times \{0\}^{n-\did}$, where $\Lambda$ is a lattice of rank $p$ and dimension $\did$.
Furthermore, the image of the vectors $e^1,e^2\dots,e^p$ forms a lattice basis $b^1,\dots,b^p$ of~$\Lambda$.
Since $\N \subseteq \L$, the image of~$\N$ is $\N' \times \{0\}^{n-\did}$, where $\N'$ is a linear subspace of~$\R^\did$ of dimension $\did-p$.
Since $\M$ and $\N$ are orthogonal, we have that $\Lambda + \N'$ has dimension $d$.
Finally, we know that the image of~$\L^\perp$ is $\{0\}^\did \times \R^{n-\did}$.
We conclude that the image of~$\Z^p \times \R^{n-p}$ is $(\Lambda + \N') \times \R^{n-\did}$.
Let $\Lambda'$ be the orthogonal projection of~$\Lambda$ onto ${\N'}^\perp$.
Then $\Lambda'$ is a lattice of rank $p$ and dimension $\did$, and the image of~$\Z^p \times \R^{n-p}$ is $(\Lambda' + \spn(\Lambda')^\perp) \times \R^{n-\did}$ as desired.
A basis of~$\Lambda'$ can be obtained by taking the orthogonal projection of~$b^1,\dots,b^p$ onto ${\N'}^\perp$.

By eventually reordering the components of the vector $y$, and accordingly the data of the problem, we obtain that the diagonal elements of the matrix $D$ satisfy $\abs{D_{11}} \ge \cdots \ge \abs{D_{\did\did}}$.
\end{prf}

Next, we briefly discuss how Proposition~\ref{prop change of basis} simplifies in the pure integer setting and in the pure continuous setting.
In the pure integer setting we have $p=n$ in \eqref{prob main}, and Proposition~\ref{prop change of basis} implies $\did = n$.
Therefore, in \eqref{prob s main} we have no $z$ variables and the constraint $y \in \Lambda + \spn(\Lambda)^\perp$ is replaced by $y \in \Lambda$ since the set $\Lambda$ is a full rank lattice of dimension $n$.
Furthermore, in \eqref{eq containment}, the set $\proj_y \P$ is replaced by $\P$.
In the pure continuous setting we have $p=0$ in \eqref{prob main}, and Proposition~\ref{prop change of basis} implies $\did = k$.
Therefore, in \eqref{prob s main} the constraint $y \in \Lambda + \spn(\Lambda)^\perp$ is replaced by $y \in \R^\did$ since the set $\Lambda$ is a lattice of rank zero.

We remark that a change of basis similar to the one given by Proposition~\ref{prop change of basis} can be obtained through the use of eigenvalue methods like the Schur decomposition~\cite{GolVan4th}, instead of our symmetric decomposition algorithm.
These techniques have been used by Vavasis 
to obtain a related change of basis for QP (see page 282 in~\cite{Vav92i}).
Unfortunately, these methods do not yield polynomial time algorithms since symmetric integer matrices can have irrational eigenvalues.

\section{Aligned vectors}
\label{sec aligned}


In this section, we introduce the notion of aligned vectors.
Given an instance of problem \eqref{prob s main}, two vectors $y^+,y^- \in \R^\did$ are said to be \emph{aligned} if $y^+,y^- \in \B(\ctr,1) \cap (2\Lambda+\spn(\Lambda)^\perp)$, and $y^+_1 - y^-_1 \ge 1$, $\sum_{i=2}^\did (y^+_i - y^-_i)^2 \le 1/4$.
The end goal of this section is to show that, if there exist two aligned vectors,
then, for every $\epsilon \in (0,1]$, it is possible to find an $\epsilon$-approximate solution to \eqref{prob s main} by solving a number of MILPs.

We begin by showing, in Lemma~\ref{lem aligned}, how aligned vectors allow us to obtain a lower bound on the gap between maximum and minimum of a separable quadratic function evaluated on the two vectors and their midpoint.
In the proof of Lemma~\ref{lem aligned} we use the following simple lemma.
The proof is that of Lemma~3 in~\cite{Vav92i}, even though our statement is slightly stronger.

\begin{lemma}
\label{lem 1d part 1}
Let $q(\lambda) = a \lambda^2 + b \lambda + c$ be a univariate quadratic function
and let $\upp, \low \in \R$.
Let $\underline q$ and $\overline q$ be the minimum and maximum values attained by $q$ on the three points $\upp, \low, (\upp + \low)/2$.
Then $\overline q - \underline q \ge \abs{a} (\upp-\low)^2/4.$
\end{lemma}

\begin{knownproof}
\begin{prf}
By eventually replacing $q$ with $-q$, we can assume without loss of generality that $a \ge 0$.
We denote by $\lambda^\bullet$ the midpoint of the segment joining $\low$ and $\upp$, \ie $\lambda^\bullet := (\upp + \low)/2.$
Then one can readily verify that
\begin{align*}
q(\low)+q(\upp)-2q(\lambda^\bullet) = a (\upp-\low)^2/2.
\end{align*}
Since $\overline q \ge q(\low)$, $\overline q \ge q(\upp)$, and $\underline q \le q(\lambda^\bullet)$, the previous inequality implies
$\overline q - \underline q \ge a (\upp-\low)^2/4.$
\end{prf}
\end{knownproof}


\begin{lemma}
\label{lem aligned}
Let $f : \R^\did \times \R^{n-\did} \to \R$ be a quadratic function of the form $f(y,z) = y^\transp D y + c^\transp y + \zcoeff^\transp z$, where $D$ is diagonal and $D_{11}$ is the element of~$D$ with the largest absolute value.
Let $(y^+,z^+), (y^-,z^-) \in \R^\did \times \R^{n-\did}$ such that 
$y^+_1 - y^-_1 \ge 1$ and
$\sum_{i=2}^\did (y^+_i - y^-_i)^2 \le 1/4$.
Let $\underline f$ and $\overline f$ be the minimum and maximum values attained by $f$ on the three vectors $(y^+,z^+), (y^-,z^-), (y^+,z^+)/2 + (y^-,z^-)/2$.
Then $\overline f - \underline f \ge \frac 3{16} \abs{D_{11}}.$
\end{lemma}

\begin{prf}
By eventually replacing $f$ with $-f$, we can assume without loss of generality that $D_{11} \ge 0$.
Let $q : \R \to \R$ be defined by 
\begin{align*}
q(\lambda) := f\pare{(y^-,z^-) + \lambda \pare{(y^+,z^+) - (y^-,z^-)}}.
\end{align*}
Using the separability of~$f$, we obtain
\begin{align*}
q(\lambda) 
= \sum_{i=1}^\did D_{ii} \pare{y^-_i + \lambda (y^+_i - y^-_i)}^2 + O(\lambda) 
= \lambda^2 \cdot \sum_{i=1}^\did D_{ii} (y^+_i - y^-_i)^2  + O(\lambda).
\end{align*}

To conclude the proof we just need to show that 
\begin{align}
\label{eq lastthing yz}
\absL{\sum_{i=1}^\did D_{ii} (y^+_i - y^-_i)^2} \ge \frac 34 D_{11}.
\end{align}
In fact, by noting that 
$q(0) = f(y^-,z^-)$,
$q(1) = f(y^+,z^+)$, 
and $q(1/2) = f \pare{(y^+,z^+)/2 + (y^-,z^-)/2}$,
we can apply Lemma~\ref{lem 1d part 1} to $q$ and the points $0,1 \in \R$ and obtain
\begin{align*}
\overline f - \underline f 
= \overline q - \underline q 
\ge \frac 14 \absL{\sum_{i=1}^\did D_{ii} (y^+_i - y^-_i)^2}
\ge \frac 3{16} D_{11}.
\end{align*}

To prove inequality \eqref{eq lastthing yz}, we bound its left hand side as follows:
\begin{align*}
\absL{\sum_{i=1}^\did D_{ii} (y^+_i - y^-_i)^2} 
& \ge \sum_{i=1}^\did D_{ii} (y^+_i - y^-_i)^2 = \\
& = \sum_{i : D_{ii} \ge 0} D_{ii} (y^+_i - y^-_i)^2 - \sum_{i : D_{ii} < 0} -D_{ii} (y^+_i - y^-_i)^2.
\end{align*}
We can now separately bound the two nonnegative sums using 
the assumption on $D_{11}$, and the conditions
$y^+_1 - y^-_1 \ge 1$ and
$\sum_{i=2}^\did (y^+_i - y^-_i)^2 \le 1/4$.
\begin{align*}
\sum_{i : D_{ii} \ge 0} D_{ii} (y^+_i - y^-_i)^2 & \ge D_{11} (y^+_1 - y^-_1)^2 \ge D_{11}, \\
\sum_{i : D_{ii} < 0} -D_{ii} (y^+_i - y^-_i)^2 & \le D_{11} \sum_{i : D_{ii} < 0} (y^+_i - y^-_i)^2 \le D_{11} / 4.
\end{align*}
Hence inequality \eqref{eq lastthing yz} holds.
\end{prf}


We are now ready to discuss our approximation algorithm for spherical form MIQPs for which there exist two aligned vectors.
%
%
%
%
This algorithm 
is based on 
the classic technique
of
mesh partition and linear underestimators.
This natural approach consists in replacing the nonlinear objective function with a piecewise linear approximation, an idea known in the field of optimization since at least the 1950s.
%
Mesh partition and linear underestimators have proven to be a very successful technique to obtain approximation algorithms for 
several special classes of MIQP~\cite{Vav92c,Vav92i,dP16,dP18,dP19}.
In this section, for the first time we employ mesh partition and linear underestimators to MIQPs that, at the same time, have integer variables and an indefinite quadratic objective function. 
The generality of this setting poses a number of additional challenges, and the results presented in the paper so far 
provide the key to successfully apply these techniques.
In the proof, we will use the following standard lemma. 

\begin{lemma}
\label{lem 1d part 2}
Let $q(\lambda) = a \lambda^2 + b \lambda + c$ be a univariate quadratic function and let $\upp, \low \in \R$.
Let $\linearizedinlemma(\lambda)$ be the affine univariate function 
that attains the same values as $q$ at $\low,\upp$.
Then $\abs{\linearizedinlemma(\lambda) - q(\lambda)} \le \abs{a} (\upp-\low)^2/4$ for every $\lambda \in [\low,\upp]$.
\end{lemma}

\begin{knownproof}
\begin{prf}
It is simple to check that 
$\linearizedinlemma(\lambda) = a (\low+\upp) \lambda - a \low \upp + b \lambda + c$.
By eventually replacing $q$ with $-q$, we can assume without loss of generality that $a \ge 0$.
The convexity of~$q$ then implies $\linearizedinlemma(\lambda) \ge q(\lambda)$ for every $\lambda \in [\low,\upp]$.
It can be checked that
$
\linearizedinlemma(\lambda) - q(\lambda) = 
a (\lambda-\low) (\upp - \lambda).
$
This univariate quadratic function is concave, and its maximum is achieved at $(\low+\upp)/2$.
This maximum value is $a(\upp - \low)^2/4$, therefore $\linearizedinlemma(\lambda) - q(\lambda) \le a(\upp - \low)^2/4$ for every $\lambda \in \R$.
\end{prf}
\end{knownproof}




\begin{proposition}
\label{prop latest}
%
Consider \eqref{prob s main}, assume that there exist two aligned vectors, and let $k$ be the rank of the matrix $D$.
For every $\epsilon \in (0,1]$, there is an algorithm that finds an $\epsilon$-approximate solution, if it exists, 
by solving at most $\ceilL{4 \constelli_{d} \sqrt{k/(3 \epsilon)}}^k$ MILPs of the same size as \eqref{prob s main} and with $p$ integer variables.
\end{proposition}

\begin{prf}
We start by describing the approximation algorithm.
We define $\constboxes^k$ boxes in~$\R^{k}$, where $\constboxes := \ceilL{4 \constelli_\did \sqrt{k/(3 \epsilon)}}$:
\begin{align}
\label{eq boxes}
\C_{j_1,\dots,j_k} := 
\prod_{i=1}^k \pare{ \bra{\ctr_i - \constelli_\did} + \frac{2 \constelli_\did}{\constboxes} \left[j_i - 1 , j_i \right] }
\quad \forall j_1,\dots,j_k \in \{1,\dots, \constboxes\}.
\end{align}
Note that the union of these $\constboxes^k$ boxes is the polytope
\begin{align*}
\{(y_1,\dots,y_k) \in \R^{k} : \ctr_i - \constelli_\did \le y_i \le \ctr_i + \constelli_\did \ \forall i =1,\dots,k\},
\end{align*}
which contains the projection of~$\P$ onto the space defined by the first $k$ coordinates of~$y$, since $\proj_y \P \subset \B(\ctr,\constelli_\did)$ from \eqref{eq containment}.

For each box $\C = \prod_{i=1}^k [\low_i,\upp_i]$ among those defined in \eqref{eq boxes}, we construct the affine 
functions
$\linearized_i : \R \to \R$ that attain the same values as $D_{ii} y_i^2$
at $\low_i,\upp_i$, for $i=1,\dots,k$:
\begin{align*}
\linearized_i(y_i) := D_{ii}(\low_i+\upp_i) y_i 
- D_{ii} \low_i \upp_i \qquad \forall i = 1,\dots, k.
\end{align*}
We define $\gamma := \abs{D_{11}}.$
Then we define the affine function $\linearized: \R^{k} \to \R$ given by
\begin{align}
\label{eq underestimator}
\linearized(y_1,\dots,y_k) := \sum_{i=1}^k \linearized_i(y_i) - \frac{\gamma \constelli_\did^2}{\constboxes^2} \abs{\{i \in \{1,\dots,k\} : D_{ii} > 0\}}.
\end{align}
We solve
the MILP obtained from \eqref{prob s main} by substituting $y^\transp D y$ with $\linearized(y_1,\dots,y_k)$ and adding the constraint $(y_1,\dots,y_k) \in \C$:
\begin{align}
\label{prob MILP on box}
\begin{split}
\min & \quad \linearized(y_1,\dots,y_k) 
+ c^\transp y
+ \zcoeff^\transp z \\
\st & \quad (y,z) \in \P  \\
& \quad (y_1,\dots,y_k) \in \C \\
& \quad y \in \Lambda + \spn(\Lambda)^\perp, \ z \in \R^{n-\did}.
\end{split}
\end{align}
To see that \eqref{prob MILP on box} is indeed a MILP, one just needs to perform a change of basis that maps $\Lambda$ to $\Z^p \times \{0\}^{d-p}$ and $\spn(\Lambda)^\perp$ to $\{0\}^{p} \times \R^{d-p}$.



The approximation algorithm returns the best solution $(y^\diamond, z^\diamond)$ among all the (at most) $\constboxes^k$ optimal solutions just obtained of the MILPs \eqref{prob MILP on box}.
If all the MILPs \eqref{prob MILP on box} are infeasible, the algorithm returns that \eqref{prob s main} is infeasible.
This concludes the description of the algorithm.

\smallskip

%
%
Next, we show that 
$(y^\diamond,z^\diamond)$ 
is an $\epsilon$-approximate solution to 
\eqref{prob s main}.
To simplify the notation, in this proof we denote the objective function of \eqref{prob s main} by 
\begin{align*}
f(y, z) := y^\transp D y + c^\transp y + \zcoeff^\transp z = 
\sum_{i=1}^k D_{ii} y_i^2
+ c^\transp y + \zcoeff^\transp z.
\end{align*}

In order to show that $(y^\diamond,z^\diamond)$ is an $\epsilon$-approximate solution, we derive two bounds:
(i)
an upper bound on $f(y^\diamond,z^\diamond) - f(y^*,z^*)$, where $(y^*, z^*)$ is an optimal solution to \eqref{prob s main}, and
(ii)
a lower bound on $f_{\max} - f(y^*,z^*)$, where $f_{\max}$ is the maximum value of~$f(y,z)$ on the feasible region of \eqref{prob s main}.
Note that both bounds will depend linearly on 
$\gamma$. 
This dependence 
is what allows us to solve a number of MILPs that is independent on $\gamma$.

\smallskip
\noindent
\ul{An upper bound on $f(y^\diamond,z^\diamond) - f(y^*,z^*)$.} \space
Let $\C \subset \R^k$ be a box constructed in \eqref{eq boxes}, say 
$\C = \prod_{i=1}^k [\low_i,\upp_i]$.
For each $i = 1,\dots, k,$ we apply Lemma~\ref{lem 1d part 2} to each univariate quadratic function
$
D_{ii} y_i^2
$ and points $\low_i, \upp_i$.
Since $\upp_i - \low_i = 2 \constelli_\did / \constboxes$ and $\abs{D_{ii}} \le \gamma$ for $i=1,\dots,k$, we obtain that, for every $(y_1,\dots,y_k) \in \C$,
\begin{align*}
\linearized_i(y_i) - \gamma \constelli_\did^2 /\constboxes^2 & \le D_{ii} y_i^2 \le \linearized_i(y_i)
&& \text{if } D_{ii} > 0\\
\linearized_i(y_i) & \le D_{ii} y_i^2 \le \linearized_i(y_i) + \gamma \constelli_\did^2 /\constboxes^2
&& \text{if } D_{ii} < 0.
\end{align*}
We sum up all these inequalities for $i=1,\dots,k$ and obtain that for every \mbox{$(y_1,\dots,y_k) \in \C$},

\begin{align}
\label{eq claim sub}
\linearized(y_1,\dots,y_k) \le 
\sum_{i=1}^k D_{ii} y_i^2
\le \linearized(y_1,\dots,y_k) + \gamma k \constelli_\did^2 / \constboxes^2.
\end{align}
Let $\C^\diamond \subset \R^k$ be the box constructed in \eqref{eq boxes} that yields the solution $(y^\diamond,z^\diamond)$ and let $\linearized^\diamond$ be the corresponding affine function defined in \eqref{eq underestimator}.
Let $\C^* \subset \R^k$ be a box constructed in \eqref{eq boxes} such that $(y^*,z^*) \in \C^*$ and let $\linearized^*$ be the corresponding affine function.
We have
\begin{align}
\begin{split}
\label{eq claim good}
f(y^\diamond,z^\diamond) 
& \le \linearized^\diamond(y_1^\diamond,\dots,y_k^\diamond) 
+ c^\transp y^\diamond 
+ \zcoeff^\transp z^\diamond + \gamma k \constelli_\did^2 / \constboxes^2 \\
& \le \linearized^*(y_1^*,\dots,y_k^*) 
+ c^\transp y^*
+ \zcoeff^\transp z^* + \gamma k \constelli_\did^2 / \constboxes^2 \\
& \le f(y^*,z^*) + \gamma k \constelli_\did^2 / \constboxes^2.
\end{split}
\end{align}
The first inequality follows by applying the right inequality in \eqref{eq claim sub} to $\C^\diamond$ and $y^\diamond$.
The second inequality holds by definition of~$(y^\diamond,z^\diamond)$.
The third inequality follows by applying the left inequality in \eqref{eq claim sub} to $\C^*$ and $y^*$.

\smallskip
\noindent
\ul{A lower bound on $f_{\max} - f(y^*, z^*)$.} \space
By assumption, there exist two aligned vectors $y^+,y^-$ 
for \eqref{prob s main}.
From \eqref{eq containment} we have $\B(\ctr,1) \subset \proj_y \P$, thus there exist $z^+,z^- \in \R^{n-\did}$ such that 
the vectors $(y^+,z^+), (y^-,z^-) \in \R^\did \times \R^{n-\did}$ are in~$\P$.
We define the midpoint of the segment joining $(y^+,z^+)$ and $(y^-,z^-)$ as $(y^\circ,z^\circ) :=  (y^+,z^+)/2 + (y^-,z^-)/2.$
By convexity, the vector $(y^\circ,z^\circ)$ is in~$\P$.
Moreover, as both vectors $y^+ /2$, $y^- /2$ are in~$\Lambda+\spn(\Lambda)^\perp$, so is their sum $y^\circ$.
Let $\underline f$ and $\overline f$ be the minimum and maximum values attained by $f$ on the three vectors $(y^+,z^+)$, $(y^-,z^-)$, $(y^\circ,z^\circ)$.
Then, by Lemma~\ref{lem aligned}, 
$
\overline f - \underline f \ge \frac 3{16} \abs{D_{11}} = \frac 3{16} \gamma.
$
Since all three vectors are feasible to \eqref{prob s main}, we conclude that 
\begin{align}
\label{eq claim bad}
f_{\max} - f(y^*, z^*) \ge \frac 3{16} \gamma.
\end{align}

\smallskip

We are now ready to show that $(y^\diamond,z^\diamond)$ is an $\epsilon$-approximate solution to \eqref{prob s main}.
We have 
\begin{align*}
\frac{f(y^\diamond,z^\diamond) - f(y^*, z^*)}{f_{\max} - f(y^*, z^*)}
& \le \frac{\cancel{\gamma} k \constelli_\did^2}{\constboxes^2} \cdot \frac{16}{3 \cancel{\gamma}} 
= \frac{16}{3} \frac{k \constelli_\did^2}{\constboxes^2} 
\le \epsilon.
\end{align*}
In the first inequality we used \eqref{eq claim good} and \eqref{eq claim bad}, and the  last inequality holds by the definition of~$\constboxes$
given at the beginning of the proof.
\end{prf}

In particular, note that the number of MILPs solved in Proposition~\ref{prop latest} is polynomial in~$1/\epsilon$ if $k$ and $d$ are fixed.
Due to Proposition~\ref{prop change of basis}, this is indeed the case if both $k$ and $p$ are fixed in the original \eqref{prob main}.

\section{Flatness and decomposition of spherical form MIQP} 
\label{sec flatness}


In Section~\ref{sec aligned} we saw that, if a spherical form MIQP has two aligned vectors, then we can find an $\epsilon$-approximate solution.
But what if there are no aligned vectors?
In this section, we show that in this case we can decompose the problem in a number of MIQPs with fewer integer variables.
%
This result will play a crucial role in our approximation algorithm for MIQP.
Before stating our theorem, we recall the concepts of width and of reduced basis.

Let $\S \subseteq \R^\did$ be a bounded closed convex set.
Given a vector $\direction \in \R^\did$, we define the \emph{width of~$\S$ along $\direction$} to be
\begin{align*}
\width_\direction (\S) := \max \{ \direction^\transp y : y \in \S \} - \min \{ \direction^\transp y : y \in \S \}.
\end{align*}
Let $\Lambda$ be a lattice of rank $p$ and dimension $\did$, and let $b^1, \dots , b^p \in \R^\did$ be a lattice basis of~$\Lambda$.
Consider now a vector $\direction \in \R^\did$ that satisfies $\direction^\transp b^i \in \Z$ for every $i =1,\dots,p$.
Then $\direction^\transp y$ is an integer for every $y \in \Lambda$ since $y$ can be written as an integer linear combination of the $b^i$.
It follows that $\width_\direction(\S)$ is an upper bound on the number of hyperplanes orthogonal to $\direction$ that contain points in~$\S \cap \Lambda$.

Next, we recall the notion of reduced basis.
Let $\Lambda$ be a lattice of rank $p$ and dimension $\did$, and let $b^1, \dots , b^p \in \R^\did$ be a lattice basis of~$\Lambda$.
The $\did \times p$ matrix $B$ formed by taking the columns to be the basis vectors $b^i$ is called a \emph{basis matrix} of~$\Lambda$.
The \emph{determinant} of~$\Lambda$ is the volume of the fundamental parallelepiped of any basis for $\Lambda$, that is, $\det(\Lambda) := \sqrt{\det(B^\transp B)}$.

Lov\'asz introduced the notion of a reduced basis, using a Gram-Schmidt orthogonal basis as a reference. 
The \emph{Gram-Schmidt procedure} is as follows. 
Define $g^1 := b^1$ and, recursively, for $i = 2, \dots, p$, define $g^i \in \R^\did$ as the projection of~$b^i$ onto the orthogonal complement of the linear space spanned by $b^1, \dots, b^{i-1}$. 
Formally, for $i = 2, \dots, p$, $g^i$ is defined by
\begin{align}
\label{eq 9.3}
g^i &:= b^i - \sum_{j=1}^{i-1} \mu_{ij} g^j,
\qquad \text{where }
\mu_{ij} := \frac{(b^i)^\transp g^j}{\norm{g^j}^2} & \qquad 
\forall j=1,\dots,i-1.
\end{align}
By construction, the Gram-Schmidt basis $g^1 , \dots , g^p$ is an orthogonal basis of~$\spn(\Lambda)$ with the property that, for $i = 1,\dots,p$, the linear spaces spanned by $b^1,\dots, b^i$ and by $g^1,\dots, g^i$ coincide.
Moreover, we have $\norm{b^i} \ge \norm{g^i}$ for $i = 1,\dots, p$,
and $\norm{g^1} \cdots \norm{g^p} = \det(\Lambda)$.

A basis $r^1, \dots,r^p$ of the lattice $\Lambda$ is said to be \emph{reduced} if it satisfies the following two conditions
\begin{align*}
\abs{\mu_{ij}}  \le \frac 12 & \qquad \text{for } 1 \le j < i \le p \\ 
\norm{g^i + \mu_{i,i-1}g^{i-1}}^2 \ge \frac 34 \norm{g^{i-1}}^2 & \qquad \text{for } 2 \le i \le p, 
\end{align*}
where $g^1, \dots , g^p$ is the 
output of the Gram-Schmidt procedure when applied to $r^1, \dots , r^p$.
Lov\'asz' celebrated \emph{basis reduction algorithm} yields a reduced basis, and it runs in polynomial time in the size of the original basis.
If a basis $r^1, \dots , r^p$ of~$\Lambda$ is reduced, then it
is ``nearly orthogonal'', in the sense that it satisfies 
\begin{align}
\label{eq ort}
\norm{r^1} \cdots \norm{r^p} \le \constlen \det(\Lambda).
\end{align}
See for example~\cite{ConCorZamBook} for more details on lattices and reduced basis, or~\cite{GalBook} for an exposition that does not consider only full rank lattices.

In order to show our decomposition result for spherical form MIQP,
we first prove the following Lenstra-type proposition. 

\begin{proposition}
\label{prop Lenstra}
Let $a \in \Q^\did$, $\delta \in \Q$ with $\delta \ge 0$, and let $\Lambda$ be a lattice of rank $p$ and dimension $\did$ with basis matrix $\basmatB \in \Q^{\did \times p}$.
There is a polynomial time algorithm which either finds a vector $\bar y \in \B(\ctr,\delta) \cap (\Lambda+\spn(\Lambda)^\perp)$, or finds a vector $\direction \in \spn(\Lambda)$ with $\direction^\transp \basmatB$ integer such that $\width_\direction(\B(\ctr,\delta)) \le p \constlen$.
\end{proposition}

\begin{prfh}
If $p=0$, then the algorithm simply returns $\bar y = a$, thus we now assume $p \ge 1$.
The basis reduction algorithm gives in polynomial time a reduced basis $r^1, \dots , r^p \in \Q^\did$ of the lattice $\Lambda$.
Let $\hat r^1, \dots \hat r^p \in \Q^\did$ be obtained by reordering $r^1, \dots r^p$ so that the vector in the last position has maximum norm, 
and denote by $\hat R \in \Q^{\did \times p}$ the corresponding basis matrix.
Since $\basmatB$ and $\hat R$ are basis matrices of the same lattice $\Lambda$, it is well known that we can find in polynomial time a $p \times p$ unimodular matrix $U$ such that $\basmatB = \hat RU$.

Let $\ctr_\Lambda := \proj_{\spn(\Lambda)} \ctr \in \Q^\did$, let $\lambda \in \Q^p$ be such that $\hat R \lambda = a_\Lambda$, and define $y_\Lambda := \hat R \round \lambda \in \Q^\did$, where $\round \lambda = (\round{\lambda_1},\dots,\round{\lambda_p})$ and $\round{\lambda_i}$ denotes a nearest integer to $\lambda_i$.
Clearly, $y_\Lambda \in \Lambda$.
Consider first the case 
$y_\Lambda \in \proj_{\spn(\Lambda)}(\B(\ctr,\delta))$.
This implies that the vector $\bar y := (\ctr + \spn(\Lambda)) \cap (y_\Lambda +\spn(\Lambda)^\perp)) \in \Q^\did$ is in~$\B(\ctr,\delta)$.
Since $y_\Lambda \in \Lambda$, we have that $\bar y \in \Lambda+\spn(\Lambda)^\perp$.
Therefore, in this case we are done. 
Hence, in the remainder of the proof we consider the case $y_\Lambda \notin \proj_{\spn(\Lambda)}(\B(\ctr,\delta))$.

Since $\basmatB$ is a $\did \times p$ matrix of rank $p$, the matrix $\basmatB^\transp \basmatB$ is an invertible $p \times p$ symmetric matrix, thus we can define the $p \times \did$ matrix $\basmatB^\dagger := (\basmatB^\transp \basmatB)^{-1} \basmatB^\transp$.
The matrix $\basmatB^\dagger$ is a left inverse of~$\basmatB$, \ie $\basmatB^\dagger \basmatB$ is the identity matrix $\id_p$.
Let $\lastrowU \in \Z^{1 \times p}$ be the last row of~$U$, and define the vector $\direction := (\lastrowU\basmatB^\dagger)^\transp \in \Q^\did$. 
We have that $\direction \in \spn(\Lambda)$, since for every vector $t \in (\spn(\Lambda))^\perp$ we have
$\direction^\transp t = \lastrowU \basmatB^\dagger t = \lastrowU (\basmatB^\transp \basmatB)^{-1} \basmatB^\transp t = 0,$
since each column of~$\basmatB$ lies in~$\spn(\Lambda)$.
Moreover, the vector $\direction^\transp \basmatB$ is integer since 
$\direction^\transp \basmatB = \lastrowU\basmatB^\dagger \basmatB = \lastrowU \id_p = \lastrowU.$
Hence, to complete the proof, we only need to show $\width_\direction(\B(\ctr,\delta)) \le p \constlen$.

The assumption $y_\Lambda \notin \proj_{\spn(\Lambda)}(\B(\ctr,\delta))$ is equivalent to $\norm{y_\Lambda - \ctr_\Lambda} > \delta$. 
Since $y_\Lambda = \hat R \round \lambda$ and $\ctr_\Lambda = \hat R \lambda$, we have
\begin{align*}
\norm{y_\Lambda - \ctr_\Lambda} 
& = \norm{\hat R(\round \lambda - \lambda)} = 
\normL{\sum_{i=1}^p (\round{\lambda_i} - \lambda_i) \hat r^i} \\
& \le \sum_{i=1}^p \absL{\round{\lambda_i} - \lambda_i} \ \norm{\hat r^i} \le p \norm{\hat r^p}/2.
\end{align*}
We obtain that $\norm{\hat r^p} > 2 \delta/p$.
Consider the Gram-Schmidt orthogonal basis $\hat g^1, \dots , \hat g^p \in \Q^\did$ obtained from $\hat r^1,\dots,\hat r^p$. 
Using \eqref{eq ort} we have 
\begin{align*}
\norm{\hat r^1} \cdots \norm{\hat r^p} = \norm{r^1} \cdots \norm{r^p} \le \constlen \det(\Lambda) = \constlen \norm{\hat g^1} \cdots\norm{\hat g^p}.
\end{align*}
Moreover, as $\norm{\hat r^i} \ge \norm{\hat g^i}$ for $i = 1,\dots,p-1$, it follows that $\norm{\hat r^p} \le \constlen \norm{\hat g^p}$. 
Since $\norm{\hat r^p} > 2 \delta / p$, we obtain
\begin{align}
\label{9.8mi}
\norm{\hat g^p} > \frac{2 \delta}{p\constlen}.
\end{align}
We define the $p \times \did$ matrix $\hat R^\dagger := (\hat R^\transp \hat R)^{-1} \hat R^\transp$, which is a left inverse of~$\hat R$.
Using $\basmatB = \hat R U$, we obtain the relation
\begin{align*}
\basmatB^\dagger 
& = (\basmatB^\transp \basmatB)^{-1} \basmatB^\transp
= (U^\transp \hat R^\transp \hat R U)^{-1} U^\transp \hat R^\transp \\
& = U^{-1} (\hat R^\transp \hat R)^{-1} U^{-\transp} U^\transp \hat R^\transp 
= U^{-1} (\hat R^\transp \hat R)^{-1} \hat R^\transp
= U^{-1} \hat R^\dagger.
\end{align*}
It is simple to check that $\width_\direction(\B(\ctr,\delta)) = 2 \delta\norm{\direction}$.
If we denote by $\lastrowR \in \Q^{1 \times \did}$ the last row of~$\hat R^\dagger$, we have
\begin{align}
\label{eq wdmi}
\width_\direction(\B(\ctr,\delta)) 
= 2 \delta\norm{\direction} 
= 2 \delta\norm{(\lastrowU\basmatB^\dagger)^\transp}
= 2 \delta\norm{(\lastrowU U^{-1}\hat R^\dagger)^\transp}
= 2 \delta\norm{\lastrowR^\transp},
\end{align}
where the last equality holds since $\lastrowU$ is the last row of~$U$. 

We now show that $\lastrowR^\transp = \hat g^p/\norm{\hat g^p}^2$.
First, note that both $\lastrowR^\transp$ and $\hat g^p$ live in~$\spn(\Lambda)$.
For $\hat g^p$ this follows from the fact that $\hat g^1,\dots,\hat g^p$ is a basis of~$\spn(\Lambda)$.
For $\lastrowR^\transp$, it can be seen because this vector is orthogonal to each vector $t \in (\spn(\Lambda))^\perp$ as $\lastrowR$ is the last row of~$\hat R^\dagger$ and we have $\hat R^\dagger t = (\hat R^\transp \hat R)^{-1} \hat R^\transp t = 0$, since each column of~$\hat R$ lies in~$\spn(\Lambda)$.
Since $\hat g^p$ is orthogonal to $\hat g^1,\dots, \hat g^{p-1}$, it follows from \eqref{eq 9.3} that $(\hat g^p)^\transp \hat r^i = 0$ for $i = 1,\dots,p-1$ and $(\hat g^p)^\transp \hat r^p = \norm{\hat g^p}^2$. 
Since $\lastrowR$ is the last row of~$\hat R^\dagger$, we have $\lastrowR \hat r^i = 0$ for $i = 1,\dots,p-1$ and $\lastrowR \hat r^p = 1$.
This concludes the proof that $\lastrowR^\transp = \hat g^p/\norm{\hat g^p}^2$.

Thus, by \eqref{eq wdmi} and \eqref{9.8mi},
\begin{equation*}
\width_\direction(\B(\ctr,\delta)) 
= 2 \delta\norm{\lastrowR^\transp}
= \frac{2 \delta}{\norm{\hat g^p}} 
\le p \constlen. \qedhere
\end{equation*}
\end{prfh}

We are now ready to give our decomposition result.

\begin{proposition}
\label{prop aligned}
There is a polynomial time algorithm which either 
finds two aligned vectors
for \eqref{prob s main},
or finds a vector $\direction \in \spn(\Lambda)$ with $\direction^\transp B$ integer such that 
$\width_\direction(\P) \le \constelli_\did \constwidth_p$, where 
$\constwidth_p := 14 p \constlen$.
\end{proposition}


\begin{prf}
Let $\ctr^+ := \ctr + \frac 34 e^1 \in \Q^\did$,
where $e^1$ denotes the first vector of the standard basis of~$\R^\did$.
It is simple to verify that
\begin{align}
\label{eq claim cont mi}
\B(\ctr^+, 1/4) \subseteq \B(\ctr,1) \subseteq \B(\ctr^+, 7/4).
\end{align}

Denote by $B \in \Q^{\did \times p}$ the given basis matrix of the lattice $\Lambda$.
We apply Proposition~\ref{prop Lenstra} to $\B(\ctr^+, 1/4)$ and the lattice $2 \Lambda$ with basis matrix $2\basmatB$.
Consider first the case where Proposition~\ref{prop Lenstra} finds a vector $\direction \in \spn(\Lambda)$ with $\direction^\transp (2\basmatB)$ integer such that $\width_\direction(\B(\ctr^+,1/4)) \le p \constlen$.
We then set $\direction' := 2\direction$ and note that $\direction' \in \spn(\Lambda)$ with ${\direction'}^\transp B$ integer.
Furthermore, it follows from \eqref{eq claim cont mi} that
\begin{align*}
\width_{\direction'}(\B(\ctr,1)) 
&= 2 \width_\direction(\B(\ctr,1)) 
\le 2 \width_\direction(\B(\ctr^+, 7/4)) \\
&= 14 \width_\direction(\B(\ctr^+,1/4)) 
\le 14 p \constlen = \constwidth_p.
\end{align*}
Using \eqref{eq containment} we obtain
\begin{align*}
\width_{\direction'}(\P)
& = \width_{\direction'}(\proj_y \P)) 
\le \width_{\direction'}(\B(\ctr,\constelli_\did)) \\
& \le \constelli_\did \width_{\direction'}(\B(\ctr,1)) 
\le \constelli_\did \constwidth_p.
\end{align*}
Hence the statement of the proposition holds. 
Therefore, we now assume that Proposition~\ref{prop Lenstra} finds 
a vector $y^+ \in \B(\ctr^+,1/4) \cap (2\Lambda+\spn(\Lambda)^\perp)$.
Clearly, \eqref{eq claim cont mi} implies that $y^+ \in \B(\ctr,1)$.

Next, we define $\ctr^- := \ctr - \frac 34 e^1 \in \Q^\did$, and we apply Proposition~\ref{prop Lenstra} to $\B(\ctr^-, 1/4)$ and the lattice $2 \Lambda$ with basis matrix $2\basmatB$.
Symmetrically, we can assume that Proposition~\ref{prop Lenstra} finds 
a vector $y^- \in 2\Lambda+\spn(\Lambda)^\perp$ that is in~$\B(\ctr^-,1/4)$ and therefore in~$\B(\ctr,1)$.

To conclude the proof, we only need to show that
the vectors $y^+,y^-$ are aligned for \eqref{prob s main}.
Since $y^+ \in \B(\ctr^+,1/4)$ and $y^- \in \B(\ctr^-,1/4)$, we obtain
$y^+_1 - y^-_1 \ge (\ctr_1 + 1/2) - (\ctr_1 - 1/2) = 1.$
For a vector $y \in \R^\did$ we denote by $y_{-1}$ the vector in~$\R^{\did-1}$ obtained by deleting the first component from $y$.
Using the triangle inequality and the fact that $a^+_{-1} = a^-_{-1} = a_{-1}$, we obtain
\begin{align*}
\sum_{i=2}^\did (y^+_i - y^-_i)^2 
& = \norm{y^+_{-1} - y^-_{-1}}^2 \le (\norm{y^+_{-1} - a_{-1}} + \norm{y^-_{-1} - a_{-1}})^2 \\
& = (\norm{y^+_{-1} - a^+_{-1}} + \norm{y^-_{-1} - a^-_{-1}})^2 
\le (\norm{y^+ - a^+} + \norm{y^- - a^-})^2 \\
& \le (1/4 + 1/4)^2 = 1/4.
\end{align*}
Hence $y^+,y^-$ are aligned for \eqref{prob s main}.
\end{prf}

\section{Approximation algorithm}
\label{sec main algorithm}

In this section, we present our approximation algorithm for \eqref{prob main} and we prove Theorem~\ref{th main}.
First, we present two lemmas that allow us to reduce the number of variables in MIQPs with polyhedra that are not full-dimensional. 
The arguments are direct extensions of those for pure integer MILPs (see, e.g., \cite{ConCorZamBook}).
Proofs are given for completeness.


\begin{lemma}
\label{lem lower dim hyperplane}
Let $a \in \Q^n \setminus \{0\}$, $\beta \in \Q$, $p \in \{0,\dots,n\}$.
There is a polynomial time algorithm that determines whether the set $\S:=\{x \in \Z^p \times \R^{n-p} : a^\transp x = \beta \}$ is empty or not.
If $\S \neq \emptyset$, the algorithm finds a vector $\bar x \in \Q^n$ and a matrix $M \in \Q^{n\times(n-1)}$ such that 
\begin{align*}
\S = \{\bar x + My : y \in \Z^p \times \R^{n-p-1}\} &\qquad \text{if $a_i \neq 0$ for some $i \in \{p+1, \dots,n\}$} \\
\S = \{\bar x + My : y \in \Z^{p-1} \times \R^{n-p}\} &\qquad \text{if $a_i = 0$ for all $i \in \{p+1, \dots,n\}$.}
\end{align*}
\end{lemma}

\begin{prfh}
First, consider the case $a_i \neq 0$ for some $i \in \{p+1, \dots,n\}$.
We can then rewrite $a^\transp x = \beta$ in the form $x_i = (\beta - \sum_{j \in \{1,\dots,n\} \setminus \{i\}} a_j x_j)/a_i$.
Since $x_i$ is a continuous variable, we obtain that $\S$ is nonempty.
We define the vector $\bar x \in \Q^n$ with entry $\bar x_i:=\beta/a_i$ and all other entries zero.
We also define the matrix $M \in \Q^{n \times n-1}$ obtained from the $n \times n$ identity matrix by replacing the $i$th row with the horizontal vector $-a^\transp /a_i$ and deleting column $i$.
With these definitions of $\bar x$ and $M$, we obtain
\begin{align*}
\S = \{\bar x + My : y \in \Z^p \times \R^{n-p-1}\}.
\end{align*}

Next, consider the case $a_i = 0$ for all $i \in \{p+1, \dots,n\}$.
Possibly by multiplying the equation $a^\transp x = \beta$ by the least common multiple of the denominators of the entries of $a$, we may assume that $a$ is an integral vector.
Possibly by dividing the equation $a^\transp x = \beta$ by the greatest common divisor of the entries of $a$, we may assume that $a$ has relatively prime entries.
If $\beta \notin \Z$, then $\S$ is empty and we are done.
Thus, we now assume $\beta \in \Z$.
Since $a_1,\dots,a_p$ are relatively prime, by Corollary 1.9 in~\cite{ConCorZamBook}, the equation $\sum_{j=1}^p a_j x_j = \beta$ has an integral solution $\tilde x \in \Z^p$, thus $\S$ is nonempty.
Furthermore, there exists a unimodular matrix $U \in \Z^{p \times p}$ such that 
${\tilde a}^\transp U = {e_1}^\transp$, where $\tilde a$ is the vector of the first $p$ coordinates of $a$, and $e_1$ denotes the first unit vector in $\R^p$. 
From the proof of Corollary 1.9 in~~\cite{ConCorZamBook}, both $\tilde x$ and $U$ can be computed in polynomial time.
If we define the matrix $N \in \Z^{p \times (p-1)}$ formed by the last $p-1$ columns of $U$, we have 
\begin{align*}
\left\{x \in \Z^p : \sum_{j=1}^p a_j x_j= \beta\right\} = \{\tilde x + Ny : y \in \Z^{p-1}\}.
\end{align*}
We define the vector $\bar x \in \Z^n$ by $\bar x_j := \tilde x_j$ for $j \in \{1,\dots,p\}$ and $\bar x_j := 0$ for $j \in \{p+1,\dots,n\}$.
We also define the matrix $M \in \Q^{n \times n-1}$ with block corresponding to the first $p$ rows and $p-1$ columns being equal to $N$, block corresponding to the last $n-p$ rows and $n-p$ columns being equal to the identity matrix $I_{n-p}$, and remaining entries zero.
Since $a_i = 0$ for all $i \in \{p+1, \dots,n\}$, we conclude 
\begin{align*}
\S = \{\bar x + My : y \in \Z^{p-1} \times \R^{n-p}\}. \qedhere
\end{align*}
\end{prfh}

\begin{lemma}
\label{lem lower dim MIQP}
Consider an instance of \eqref{prob main} with a nonempty feasible region.
There is a polynomial time algorithm that determines whether $\{x \in \R^n : Wx \le w \}$ is full-dimensional.
If not, it rewrites the instance as an instance of \eqref{prob main} with one fewer variable.
\end{lemma}

\begin{prfh}
It is well-known \cite{ConCorZamBook} that there is a polynomial time algorithm that determines whether $\{x \in \R^n : Wx \le w \}$ is full-dimensional, and if not, finds a rational hyperplane $\{x \in \R^n : a^\transp x = \beta \}$ that contains it.
In the latter case, we let $\bar x \in \Q^n$ and $M \in \Q^{n\times(n-1)}$ from Lemma~\ref{lem lower dim hyperplane}, and we define $H' := M^\transp H M$, $h' := 2M^\transp H^\transp \bar x + M^\transp h$, $c := \bar x^\transp H \bar x + h^\transp \bar x$, $W' := WM$, $w' := w - W \bar x$.
By Lemma~\ref{lem lower dim hyperplane}, our instance of \eqref{prob main} can be rewritten as
\begin{align*}
\min & \quad x^\transp H' x + h'^\transp x + c \\
\st & \quad W'x \le w' \\
& \quad x \in \Lambda,
\end{align*}
where 
\begin{align*}
\Lambda := 
\begin{cases}
\Z^p \times \R^{n-p-1} &\qquad \text{if $a_i \neq 0$ for some $i \in \{p+1, \dots,n\}$} \\
\Z^{p-1} \times \R^{n-p} &\qquad \text{if $a_i = 0$ for all $i \in \{p+1, \dots,n\}$.\qedhere} 
\end{cases}
\end{align*}
\end{prfh}

\subsection{Description of the approximation algorithm}
\label{sec main algorithm description}

We are now in a position to present our approximation algorithm for \eqref{prob main}.
We will make use of Proposition~\ref{prop change of basis}, 
Proposition~\ref{prop latest}, 
Proposition~\ref{prop aligned},
and Lemma~\ref{lem lower dim MIQP}. 

The input of the algorithm consists of an instance of a bounded MIQP.
Theorem~4 in~\cite{dPDeyMol17} implies that, if there is an optimal solution, there is one of size bounded by an integer $\psi$, which is polynomial in the size of the input MIQP. \footnote{Even though Theorem~4 in~\cite{dPDeyMol17} does not give $\psi$ explicitly, a formula for $\psi$, as a function of the size of the MIQP instance, can be derived from its proof.}
Therefore, we obtain an equivalent MIQP instance by restricting each variable to the segment $[-2^\psi, 2^\psi]$.
The size of the latter instance is polynomial in the size of the former.
Furthermore, it is simple to check that an $\epsilon$-approximate solution to the latter instance is also an $\epsilon$-approximate solution to the former, for every $\epsilon \in [0,1]$.
Therefore, we can now assume that our input MIQP has a bounded feasible region.

We initialize the set $\instances$ of MIQP instances to be solved as a set containing only our input MIQP, and the set of possible approximate solutions as $\solutions := \emptyset$.
Throughout the algorithm, each instance in~$\instances$ will be our input MIQP with a number of additional linear equality constraints.
On the other hand, the set $\solutions$ will contain a number of feasible solutions to the input MIQP.

\setcounter{step}{0} 

\step{Output, feasibility, full-dimensionality, and linear case}
\label{step prelim}

\noindent
\textbf{Output.} 
If $\instances = \emptyset$, then we return the solution in~$\solutions$ with the minimum objective function value if $\solutions \neq \emptyset$, and we return ``infeasible'' if $\solutions = \emptyset$.
Otherwise $\instances \neq \emptyset$, we choose a MIQP instance in~$\instances$ and we remove it from $\instances$.
Without loss of generality, the chosen MIQP instance is \eqref{prob main}.

\noindent
\textbf{Feasibility.} 
Using Lenstra's algorithm~\cite{Len83}, we check if the feasible region $\{x \in \Z^p \times \R^{n-p} : Wx \le w\}$ of \eqref{prob main} is the emptyset.
If so, we go back to Step~\ref{step prelim}.
Otherwise, \eqref{prob main} is feasible and we continue. 

\noindent
\textbf{Full-dimensionality.}
We apply recursively Lemma~\ref{lem lower dim MIQP} until the polyhedron describing the feasible region is full-dimensional.
For ease of notation, we denote the obtained instance again by \eqref{prob main}, and 
we now assume that $\{x \in \R^n : Wx \le w\}$ is full-dimensional.

\noindent
\textbf{Linear case.} 
Let $k$ be the rank of the matrix $H$.
If $k=0$, \eqref{prob main} is a MILP, and we find an optimal solution using Lenstra's algorithm.
We construct the corresponding feasible solution to the input MIQP by inverting the linear transformation just performed in ``Full-dimensionality'', and we add it to $\solutions$.
Otherwise, we have $k \ge 1$ and we continue.

\step{Reduction to spherical form}
\label{step sform}

\noindent
By Proposition~\ref{prop change of basis}, we perform a change of basis that transform \eqref{prob main} in spherical form \eqref{prob s main}, where $\did$ satisfies $\did \le k+p$, the rank of the matrix $D$ is $k$, and $\constelli_\did$ in \eqref{eq containment} is the ceiling of~$2 \did^{3/2} \ceil{(5d)^{d/2}}^2$.

Let $B \in \Q^{\did \times p}$ be the obtained basis matrix of the lattice $\Lambda$.
By Proposition~\ref{prop aligned}, we either find two aligned vectors $y^+,y^-$ 
for \eqref{prob s main},
or we find a vector $\direction \in \spn(\Lambda)$ with $\direction^\transp B$ integer such that 
$\width_\direction(\P) \le \constelli_\did \constwidth_p$,
where $\constwidth_p = 14 p \constlen$.
In the first case, continue with Step~\ref{step mesh};
In the second case, go to Step~\ref{step dec}.


\step{Mesh partition and linear underestimators}
\label{step mesh}

\noindent
By Proposition~\ref{prop latest} we obtain an $\epsilon$-approximate solution $(y^\diamond,z^\diamond)$ to 
\eqref{prob s main}.
This requires solving, with Lenstra's algorithm, at most $\ceilL{4 \constelli_{d} \sqrt{k/(3 \epsilon)}}^k$ MILPs of the same size as \eqref{prob s main} and with $p$ integer variables.
%
%
%
We construct the corresponding 
$\epsilon$-approximate solution $x^\diamond$ to the \eqref{prob main} chosen at the beginning of this iteration of the algorithm by inverting the linear transformations in Step~\ref{step sform} and in Step~\ref{step prelim}, and we add it to $\solutions$.
Then, we go back to Step~\ref{step prelim}.

\step{Decomposition}
\label{step dec}

\noindent
Since $\width_\direction(\P) \le \constelli_\did \constwidth_p$,
each point $(y,z) \in \P$ with $y \in \Lambda + \spn(\Lambda)^\perp$ is contained in one of the following polytopes:
\begin{align*}
\P_t := \{(y,z) \in \P : \direction^\transp y = t\} & \qquad \text{for } t= \ceil{\mu}, \dots, \floor{\mu + \constelli_\did \constwidth_p},
\end{align*}
where $\mu := \min\{\direction^\transp y : y \in \P\}$. 

For each $t= \ceil{\mu}, \dots, \floor{\mu + \constelli_\did \constwidth_p}$, we consider the instance obtained from \eqref{prob s main} by replacing the polytope $\P$ with $\P_t$, and 
we add to $\instances$ the MIQP obtained by inverting the linear transformations in Step~\ref{step sform} and in Step~\ref{step prelim}.
Note that the instances that we just added to $\instances$ differ from the one chosen at the beginning of this iteration of the algorithm only by the additional constraint obtained from $\direction^\transp y = t$ by inverting the two linear transformations.
Finally, we go back to Step~\ref{step prelim}.

\subsection{Analysis of the algorithm} 
\label{sec main algorithm analysis}


First, we show that the algorithm described in Section~\ref{sec main algorithm description} is correct.




\begin{claim}
The algorithm in Section~\ref{sec main algorithm description} returns an $\epsilon$-approximate solution, if it exists. 
\end{claim}

\begin{prf}
Clearly, if the input MIQP is infeasible, the algorithm correctly detects it in Step~\ref{step prelim}, thus we now assume that it is feasible.
In this case, we need to show that the algorithm returns an $\epsilon$-approximate solution to the input MIQP.
To prove this, we only need to show that the algorithm eventually adds to the set $\solutions$ an $\epsilon$-approximate solution $x^\epsilon$ to the input MIQP.
In fact, the vector returned at the end of the algorithm has objective value at most that of~$x^\epsilon$, and so it is an $\epsilon$-approximate solution to the input MIQP as well.

Let $x^* \in \R^n$ be an optimal solution to the input MIQP.
Let MIQP$^*$ be an instance stored at some point in~$\instances$ that contains in the feasible region the vector $x^*$.
Among all these possible instances, we assume that MIQP$^*$, after the ``Full-dimensionality'' transformation in Step~\ref{step prelim},
has a minimal number of integer variables.
Note that MIQP$^*$ does not get decomposed in Step~\ref{step dec}.
Otherwise, the vector $x^*$ would be feasible for one of the instances generated in Step~\ref{step dec} from MIQP$^*$, which after the ``Full-dimensionality'' transformation will have fewer integer variables than MIQP$^*$.
Hence, when the algorithm selects MIQP$^*$ from $\instances$, it performs Step~\ref{step mesh} of the algorithm, and so by Proposition~\ref{prop latest} it adds to $\solutions$ a vector $x^\epsilon$ that is an $\epsilon$-approximate solution to MIQP$^*$.
Since the feasible region of MIQP$^*$ is contained in the feasible region of the input MIQP, and since the vector $x^*$ is feasible for MIQP$^*$, it is simple to check that $x^\epsilon$ is an $\epsilon$-approximate solution to the input MIQP.
\end{prf}

%
%


We complete the proof of Theorem~\ref{th main} by showing that the running time of the algorithm matches the one stated in Theorem~\ref{th main}.

\begin{claim}
The running time of the algorithm 
in Section \ref{sec main algorithm description} is polynomial in the size of the input and in~$1/\epsilon$, provided that the rank $k$ of the matrix $H$ and the number of integer variables $p$ are fixed numbers.
\end{claim}


\begin{prf}
First, we show that the algorithm performs at most 
$(\constelli_{k+p} \constwidth_p +1)^{p+1}$ iterations, which is a fixed number if both $k$ and $p$ are fixed.
Note that the number of iterations coincides with the total number of instances that are stored in~$\instances$ throughout the execution of the algorithm.
Instances are added to $\instances$ only in Step~\ref{step dec}, where 
the MIQP chosen in that iteration gets replaced in~$\instances$ with at most $\constelli_{k'+p'} \constwidth_{p'} +1$ new instances.
Here, $k'$ denotes the rank of the quadratic objective and $p'$ denotes the number of integer variables of the chosen instance after the ``Full-dimensionality'' transformation in Step~\ref{step prelim}.
In the new instances added to $\instances$, the rank of the quadratic objective is at most $k'$, and the number of integer variables is at most $p'-1$.
In particular, this implies that for every chosen instance we have $k' \le k$ and $p' \le p$.
Finally, note that Step~\ref{step dec} may be triggered only if $p' \ge 1$.
Therefore, the total number of MIQPs that are eventually stored in~$\instances$ is at most 
$\sum_{j=0}^p (\constelli_{k+p} \constwidth_p +1)^{j} \le (\constelli_{k+p} \constwidth_p +1)^{p+1}.$

It is simple to check that each instance constructed by the algorithm and each number generated has size polynomial in the size of the input MIQP. 
Thus, to conclude the proof we only need to analyze the running time of a single iteration of the algorithm.
Each MILP encountered (in Step~\ref{step prelim} and Step~\ref{step mesh}) has at most $p$ integer variables.
Since $p$ is fixed, they can be solved with Lenstra's algorithm~\cite{Len83} in time polynomial in the size of the input MIQP. 
Step~\ref{step prelim} of the algorithm can then be performed in time polynomial in the size of the input MIQP. 
By Proposition~\ref{prop change of basis} and Proposition~\ref{prop aligned}, also Step~\ref{step sform} can be performed in time polynomial in the size of the input MIQP. 
In Step~\ref{step mesh}, the algorithm solves at most $\ceilL{4 \constelli_{k+p} \sqrt{k/(3 \epsilon)}}^k$ MILPs with at most $p$ integer variables.
Since $k$ and $p$ are fixed, this number is polynomial in~$1/\epsilon$.
Therefore, Step~\ref{step mesh} of the algorithm can be performed in time polynomial in the size of the input MIQP 
and in~$1/\epsilon$.
Step~\ref{step dec} only solves one linear program to find $\mu$ and stores 
at most $\constelli_{k+p} \constwidth_p +1$ MIQPs, which is a fixed number if both $k$ and $p$ are fixed.
\end{prf}

\bigskip

\begin{small}
\noindent
\textbf{Funding: }A. Del~Pia is partially funded by  ONR grant N00014-19-1-2322. Any opinions, findings, and conclusions or recommendations expressed in this material are those of the authors and do not necessarily reflect the views of the Office of Naval Research.
\end{small}

\ifthenelse {\boolean{MPA}}
{
\bibliographystyle{spmpsci}
}
{
\bibliographystyle{plain}
}


\end{document}